\documentclass[a4paper,10pt]{amsart}%aims

\usepackage{graphicx, graphics, epsfig, color,bbm}

\newtheorem{thm}{Theorem}[section]
\newtheorem{cor}[thm]{Corollary}
\newtheorem{lem}[thm]{Lemma}
\newtheorem{prop}[thm]{Proposition}
\newtheorem{rem}[thm]{Remark}
\theoremstyle{definition}

\newcommand{\R}{\mathbb{R}}
\newcommand{\N}{\mathbb{N}}

\newcommand{\ud}{\mathrm{d}}
\newcommand{\mc}{\mathcal}
\newcommand{\be}{\begin{equation}}
\newcommand{\ee}{\end{equation}}
\newcommand{\bs}{\begin{split}}
\newcommand{\es}{\end{split}}
\newcommand{\bee}{\begin{equation*}}
\newcommand{\eee}{\end{equation*}}

\newcommand{\eps}{\varepsilon}
\newcommand{\e}{\varepsilon}

\newcommand{\vphi}{\varphi}
\numberwithin{equation}{section}

\newcommand{\kn}{\textrm{Kn}}
\newcommand{\ma}{\textrm{Ma}}
\newcommand{\re}{\textrm{Re}}

\begin{document}

%opening
\title[Incompressible Euler limit for a gas of fermions]{The incompressible Euler limit of the Boltzmann equation for a gas of fermions}
\author{Thibaut Allemand}
\address{DMA, \'Ecole Normale Sup\'erieure, 45 rue d'Ulm, 75230 Paris Cedex 05, France \\
CEA, DAM; DIF, F-91297 Arpajon, France}

\maketitle

\begin{abstract}
We are interested in the hydrodynamic limit of the Boltzmann equation for a gas of fermions in the incompressible Euler regime. We use the
relative entropy method as improved by Saint-Raymond in the classical case \cite{saintray}. Our result is analogous to the classical case
result, but the treatment is slightly complicated by the cubic nonlinearity of the collision operator. 
\end{abstract}

\section{Introduction}

The study of quantum gases has been given increasing interest in the literature over the last decade. In particular, quantum kinetic theory
is an expanding field of research. By quantum gases we mean gases made of quantum particles, that are, bosons or fermions. The first ones
aim at aggregating together to form the so-called "Bose-Einstein condensates". Conversely, fermions obey Pauli's exclusion principle, which
prevents any pair of fermions from being in the same state. 

Among the possible models for quantum gases, the Boltzmann-type models first proposed by Nordheim in 1928 \cite{nordheim_kinetic_1928} then
Uehling and Uhlenbeck in 1933 \cite{uehling_transport_1933} are very popular. Although their range of validity is not clear, they seem to
capture some aspects of the behaviour of bosonic or fermionic particles. Indeed, it has been proved, in the bosonic (and space homogeneous)
case, that under a threshold temperature, a condensate occurs in infinite time \cite{escobedo_entropy_2005, escobedo_homogeneous_2003}. On
the contrary, solutions of the fermionic Boltzmann equation satisfy a natural $L^\infty$ bound which reflects Pauli's exclusion principle and makes the Cauchy problem easier.
There is no need for renormalized solutions in this case since one can prove the existence of global weak solutions \cite{lions,
dolbeault, moi}.

 In the present work, we investigate the incompressible Euler limit of the Boltzmann equation for a gas of fermions (or Boltzmann-Fermi equation).

\subsection{The Boltzmann equation for fermions}

It reads in nondimensional form

\be
\label{BFD}
 \ma\partial _t f+v.\nabla_x f = \frac{1}{\kn} Q(f)
\ee
where $f(t,x,v)$ is the density of particles which at time $t\in\R_+$ are at the point $x\in\R^3$ of space with velocity $v\in\R^3$. The left-hand
term expresses the free transport of particles in absence of interactions, while the right-hand term takes into account the effects of
collisions between particles. The Knudsen and Mach numbers are positive constants determined by the physical situation of the gas and
defined by
\[
 \kn=\frac{\textrm{mean free path}}{\textrm{observation length scale}}
\]
\[
 \ma=\frac{\textrm{bulk velocity}}{\textrm{speed of sound}}.
\]
The Knudsen number is a measure of the rarefaction of the gas whereas the Mach number measures its compressibility.

The collision integral $Q(f)$ is given by
\be
\label{collision integral}
 Q(f)=\int_{S^2}\int_{\R^3}B(v-v_*,\omega) \left(f'f_*'(1-f)(1-f_*)-ff_*(1-f')(1-f_*') \right)\ud v_* \ud\omega
\ee
with the usual notations
\[
 f_*=f(t,x,v_*),\quad f'=f(t,x,v'),\quad f_*'=f(t,x,v_*')
\]
and where the precollisional velocities $(v',v_*')$ are deduced from the postcollisional ones by the relations
\[
 \begin{cases}
  &v'=v-(v-v_*).\omega \omega,\\
  &v_*'=v_*+(v-v_*).\omega \omega,
 \end{cases}
\]
$\omega$ being a unit norm vector. These relations express the conservation of momentum and kinetic energy at each collision. 

The collision integral \eqref{collision integral} differs from the classical one by the presence of cubic terms, due to Pauli's
exclusion principle, and thus takes into account the quantum nature of the gas. Indeed the probability for a particle with velocity $v'$ of
taking velocity $v$ after a collision is all the more penalized as other particles at the same point already have velocity $v$. This
prevents any pair of two fermions from being in the same quantum state. As a consequence the solution $f$ satisfy a natural $L^\infty$ bound
: if the initial datum $f_0$ is chosen such that
\[
 0\le f_0 \le 1,
\]
this bound is preserved by the solution $f$ for all times. This property is of great help in the study of this equation and implies
 very different features from the classical or the bosonic Boltzmann equation even if they look very similar.

The function $B(z,\omega)$, known as the collision kernel, is
measurable, a.e. positive, and depends only on $|z|$ and on the scalar product $z\cdot \omega$. It is often assumed to satisfy Grad's cutoff
assumption:
\be
\label{Grad1}
 0<B(z,\omega)\leq C_B(1+|z|)^\beta~~\textrm{a.e. on}~S^2\times\R^3
\ee
\be
\label{Grad2}
 \int_{S^2}B(z,\omega)\ud\omega\geq \frac{1}{C_B}\frac{|z|}{1+|z|}~~\textrm{a.e. on }\R^3
\ee

for some constants $C_B>0$ and $\beta \in [0,1]$. These assumptions guarantee the existence of a solution to equation \eqref{BFD} in the
whole domain $\R^3$ or in the torus $\mathbb T^3$ \cite{dolbeault,lions}. However, in order to have existence of solutions in more general
domains and to ensure that they satisfy the local conservation of mass, momentum, and energy, as well as the entropy inequality, we will
have to make more restrictive assumptions \cite{moi}, as will be explained later.

The symmetry properties of the collision operator, coming from the fact that the transformations $(v,v_*)\mapsto (v',v_*')$ and
$(v,v')\mapsto(v_*,v_*')$ have unit jacobian, imply that, at least formally,
\be
\label{symetry}
\begin{split}
 &\int_{\R^3} Q(f) \vphi \ud v = \frac{1}{4}\int_{\R^3}\int_{S^2} B \left(f'f_*'(1- f)(1-f_*)-ff_*(1- f')(1- f_*')\right)\\
&\qquad\qquad\qquad\qquad\qquad\qquad\qquad\qquad\qquad \times \left(\vphi+\vphi_*-\vphi'-\vphi_*'\right)\ud v_* \ud\omega
\end{split}
\ee
and then
\[
 \int_{\R^3}(1,v,|v|^2)Q(f)\ud v=0.
\]
As a consequence, the solution $f$ of equation \eqref{BFD} formally satisfies the conservations of mass, momentum and kinetic
energy:
\[
\begin{cases}
 &\displaystyle{\partial_t\int f \ud v + \frac{1}{\ma}\nabla_x.\left(\int vf\ud v\right) =0}\\
 &\displaystyle{\partial_t \int v f\ud v +\frac{1}{\ma}\nabla_x.\left(\int v\otimes v f \ud v\right)=0}\\
 &\displaystyle{\partial_t \int |v|^2 f\ud v +\frac{1}{\ma}\nabla_x. \left(\int v|v|^2f\ud v\right)=0}.
\end{cases}
\]
Moreover, taking $\vphi=\log\left(\frac{f}{1- f}\right)$ in \eqref{symetry} leads to the so-called H-theorem, which expresses the second
principle of thermodynamics:
\[
 \partial_t\int_{\R^3}s(f)(v)\ud v+\frac{1}{\ma}\nabla_x.\left(\int_{\R^3}vs(f)(v)\ud v\right)=-\frac{1}{\ma\kn}D(f)
\]
with
\[
 s(f)=f\log f +(1-f)\log(1-f)
\]
and
\bee
\begin{split}
 D(f)=\frac{1}{4}\int_{\R^3\times S^2}B(v-v_*,\omega)&\left(f'f_*'(1-f)(1-f_*)-ff_*(1-f')(1-f_*')\right)\\
&\cdot\log \left(\frac{f'f_*'(1-f)(1-f_*}{ff_*(1-f')(1-f_*')}\right)\ud v_* \ud\omega.
\end{split}
\eee
Notice that $D(f)$ is non-negative. The minimizers of the entropy, which are also the functions that cancel the collision operator, are
given by :

\begin{prop}
\label{planckian}
Assume that $g\in L^1(\R^3)$ is such that $Q(g)$ and $D(g)$ are well defined and satisfy the bounds $0\leq g\le 1$. Then
 \[
  Q(g)=0\iff D(g)=0 \iff g=\frac{M}{1+M}~\textrm{or}~g=\mathbbm{1}_{v\in\Lambda}
 \]
where $M$ is a maxwellian distribution, that is,
\[
 M=a e^{-\frac{|v-u|^2}{2b}}
\]
with $a\ge 0$, $b>0$, and for some subset $\Lambda\subset\R^3$. In the first case the distributions
are called Fermi-Dirac or Planckian distributions, whereas they are called degenerate Fermi-Dirac distributions in the second case.
\end{prop}
Proofs and details are to be found in  \cite{dolbeault,escobedo_entropy_2005}. The coefficiens of the Planckian distribution $a,b,u$ are
fully determined by macroscopic parameters of the fluid with distribution $g$, namely the mass $\rho_g$, bulk velocity $u_g$ and pressure
$p_g$, defined by
\[
 \rho_g(t,x) = \int_{\R^3}g(t,x,v) \ud v,\qquad  u_g(t,x)  =\frac{1}{\rho_g(t,x)} \int_{\R^3}v g(t,x,v)\ud v
\]
and
\[
 p_g(t,x,v)=\int_{\R^3}|v-u_g(t,x)|^2g(t,x,v)\ud v.
\]
This relation will be detailed in the next Section.

\subsection{The formal hydrodynamic limit}

%We will mainly rely on \cite{saintray} and . 

In the fast relaxation regime, that is, when the Knudsen number goes to 0, the collision process becomes predominant and we expect in view
of \eqref{BFD} that a local thermodynamic equilibrium is reached almost instantaneously. The distribution $f$ describing the gas is then
close to a Planckian distribution fully determined by some coefficients $a,b,u$ which in turn are caracterized by the hydrodynamic fields
$\rho_f(t,x)$, $u_f(t,x)$ and $p_f(t,x)$. The conservation laws then become up to order $O(\kn)$
\[
\label{eulercomp}
 \begin{split}
  &\ma\partial_t \rho_f + \nabla_x \cdot(\rho_f u_f)=0 \\
  &\ma \partial_t (\rho_f u_f) + \nabla_x\cdot \left(\rho_f u_f\otimes u_f + \frac{1}{3}p_fI_3\right)=0 \\
  & \ma \partial_t (\rho_f u_f^2 + p_f) + \nabla_x \cdot\left(\rho_f |u_f|^2 u_f + \frac{5}{3}p_f u_f\right)=0
 \end{split}
\]
which are the compressible Euler equations for perfect gases ($I_3$ stands for the three dimensional identity matrix). If now the Mach number also goes to $0$, that is, in incompressible
regime, the first equation becomes $\nabla_x \cdot (\rho_f u_f)=0$, which is nothing but the incompressibility constraint. The other
equations of motions are obtained by a systematic multiscale expansion, depending on another parameter: the Reynolds number, defined by
\[
 \re = \frac{\ma}{\kn}
\]
for perfect gases, which measures the viscosity of the gas. In all the sequel, we are interested in the inviscid incompressible regime, so
that we consider
\[
   \ma=\e,\qquad \kn=\e^q,\quad q>1
\]
and investigate the asymptotic $\e\to 0$. In this scaling, $\kn=\e^{1-q}\to\infty$. The Boltzmann equation for fermions now writes
\be
\label{scaled}
\eps \partial_t f +v.\nabla_x f =\frac{1}{\eps^q}Q(f).
\ee
 Moreover, the hydrodynamic fields will be assumed
to be fluctuations around a global constant equilibrium $(\rho_0,0,p_0)$, so that, denoting the fluctuations around the mass, momentum and
pressure by $\tilde\rho$, $\tilde u$ and $\tilde p$, 
\[
 \rho_f = \rho_0 +\e \tilde \rho,\quad u_f=\e \tilde u,\quad p_f=p_0+\e\tilde p.
\]
We also define the temperature $T_f$ by
\[
 p_f=\rho_f T_f
\]
and its fluctuation $\tilde T$ by $T_f=T_0+\e\tilde T$ where $T_0$ is a constant equilibrium value. Plugging these expressions into the
hydrodynamic equations \eqref{eulercomp} we get at leading order
\[
 \nabla_x \tilde u=0,\qquad \nabla_x (T_0 \tilde\rho+ \rho_0 \tilde T)=0
\]
which are known as the incompressibility and the Boussinesq relations; at next order it comes
\[
 \begin{split}
  &\partial_t \tilde\rho + \tilde u \cdot\nabla_x \tilde\rho =0 \\
 &\partial_t \tilde u + \tilde u \cdot\nabla_x\tilde u + \nabla_x (\tilde\rho\tilde T)=0.
 \end{split}
\]

The challenge is now to make rigorous this formal limit.

The rest of the paper is organized as follows: in Section \ref{sec:properties}, we recall some facts about equation \eqref{scaled}. In
Section \ref{sec:strategy} we explain our strategy which is based on the modulated entropy and we give our main result (Theorem
\ref{limite}).  In Section \ref{sec:modulated} we compute the time
derivative of the modulated entropy. In Section \ref{sec:construction} we construct an approximate solution of \eqref{scaled} the parameters
of which satisfy the incompressible Euler equation. We give some useful intermediate results and estimates in Section \ref{intermediate}. Finally 
in Section \ref{sec:end} we end the proof of the main result by controlling some flux terms.

\section{Some details on the Boltzmann equation for fermions}
\label{sec:properties}

In this Section we explain the relation between the hydrodynamic fields and the coefficients $a,u,b$ of the Planckian distribution defined
in Proposition \eqref{planckian}. We then state the existence result to be used in the sequel.

\subsection{Relations between the coefficients and the macroscopic quantities}
\label{relation}
We recall here the definition of the hydrodynamic fields: for a given distribution function $f$ the mass $\rho_f$, the
bulk velocity $u_f$ and the pressure $p_f$ are defined by
\[
 \rho_f(t,x) = \int_{\R^3}f(t,x,v) \ud v,\qquad  u_f(t,x)  =\frac{1}{\rho_f(t,x)} \int_{\R^3}v f(t,x,v)\ud v
\]
and
\[
 p_f(t,x,v)=\int_{\R^3}|v-u_f(t,x)|^2f(t,x,v)\ud v.
\]
When the thermodynamic equilibrium is reached, the density function is an equilibrium solution
\[
 f= \frac{M}{1+M},\qquad M=a e^{-\frac{|v-u|^2}{2b}},
\]
with $a,b\ge 0$. It is then clear that $u=u_f$. Let us define for comodity the functions
\[
 F_p(a)=\int_{\R^3}|v|^p\frac{a e^{-\frac{|v|^2}{2}}}{1+a e^{-\frac{|v|^2}{2}}}.
\]
They are well defined for $a\in[0,+\infty)$. Then we write
\[
 \rho_f= b^{3/2}F_0(a),\qquad p_f=b^{5/2}F_2(a).
\]
If we define now the internal energy
\[
 e_f = \frac{1}{2\rho_f}(p_f-4\rho_f|u_f|^2),
\]
we can quote the following result from \cite{arlotti}:

\begin{prop}
There exists a positive constant $l$ such that the mapping
\[
 (a,b)\in  (0,+\infty)^2 \mapsto (\rho_f, e_f)\in \mc E,
\]
where $\mc E = \{(\rho_f, e_f) ~\textrm{s.t.}~e_f>l\rho_f^{2/3}\}$, is one-to-one.
\end{prop}

From now on, we will denote by $P_0$ a global equilibrium defined by
\be
\label{P_0}
 P_0=\frac{M_0}{1+M_0},\qquad M_0=e^{a_0} e^{-\frac{|v|^2}{2b_0}}
\ee
where $a_0, b_0 >0$ and $a_0$ is such that 
\[
 F_2(e^{a_0})>2lF_0(e^{a_0})^{\frac{5}{3}}.
\]
This ensures that perturbating $a_0$ and $b_0$ will leave the corresponding perturbed mass and internal energy inside $\mc E$.

\subsection{Existence theory}

Let $\Omega$ be a subset of $\R^3$ regular enough such that the normal is well defined on the boundary. The equation
\be
\label{BBFD}
 \partial _t f +v.\nabla_x f = Q(f)
\ee
must be supplemented with an initial condition:
\be
\label{initial}
 f(0,x,v)=f_0(x,v) \qquad\forall (x,v)\in \Omega\times \R^3,
\ee
and a boundary condition we choose to be specular reflection for simplicity:
\be
\label{boundary}
 f(t,x,v)=f(t,x,R_x(v))\qquad\forall (x,v)\in \partial\Omega\times\R^3 ~~\textrm{such that}~~n(x).v<0
\ee
where $n(x)$ is the outer unit vector normal to the boundary and $R_x(v)$ is the specular reflection law
\be
\label{specular}
 R_x(v)=v-2(v.n(x))n(x).
\ee
The choice of the specular reflection as a boundary condition makes all the boundary terms vanish in the weak formulation of the equation.
It also cancels the Prandtl layer along the boundary.

The following existence result was proved in \cite{dolbeault, moi}:

\begin{thm}
Let $\Omega$ be either $\R^3$ or a regular subset of $\R^3$. Let the collision kernel $B$ be such that
\be
\label{hyp1B}
0\le B\in L^1(\R^3\times S^2),
\ee 
and let 
\be
\label{hypf00}
f_0 \in L^\infty(\Omega\times\R^3), \qquad0\leq f_0 \leq 1.
\ee
Then, the problem \eqref{BBFD}-\eqref{boundary} has a unique solution $f$ satisfying
\[
 f\in L^\infty(\R_+\times\Omega\times\R^3),\qquad 0\leq f \leq 1 ~~\textrm{a.e.}
\]
Moreover, $f$ is absolutely continuous with respect to $t$.
\end{thm}
The assumption made on the collision kernel \eqref{hyp1B} is very strong, much more than Grad's cut-off assumption
\eqref{Grad1}-\eqref{Grad2} which would be enough to ensure the existence of a solution \cite{lions}. However, it is necessary (at the time of the writing)
in order to prove that the local conservation laws are satisfied by the solution. It is the case if we assume that the collision kernel has
the symmetry
\be
\label{hypsurB}
B(w,\omega)=q(|w|,|w.\omega|),
\ee
which is physically relevant.  The following proposition was proved in \cite{moi} :
%%  using a dispersion argument which states roughly that if
%% $f$ is the solution to the transport equation
%% \[
%%  \partial_t f + v\cdot\nabla_x f=g
%% \]
%% and $\int_0^T\int_{\Omega\times\R^3}|v|^2 g\ud v\ud x\ud t \le C_0$, then $\int_0^T \int_K\int_{\R^3}|v|^3f\ud v\ud x\ud t\le C_K$ for all
%% compact subsets $K\in\R^3$.

\begin{prop}
%\label{theoremBFDexist}
Assume that the collision kernel $B$ is as in \eqref{hyp1B} and \eqref{hypsurB}. Assume moreover that the initial datum is as in \eqref{hypf00}, and
\[
 \iint_{\Omega\times\R^3}(1+|v|^3)f_0(x,v)\ud x\ud v <+\infty.
\]
Then, the solution $f$ to \eqref{BBFD}-\eqref{boundary} satisfies, in the distributional sense, the local conservation laws
\begin{equation}
\label{conslaws}
\begin{cases}
 &\displaystyle{\partial_t \int_{\R^3}f\ud v + \nabla_x.\left(\int_{\R^3}vf\ud v\right)=0}\\
 &\displaystyle{\partial_t \int_{\R^3}vf\ud v + \nabla_x.\left(\int_{\R^3}v\otimes vf\ud v\right)=0}\\
 &\displaystyle{\partial_t \int_{\R^3}|v|^2f\ud v + \nabla_x.\left(\int_{\R^3}v|v|^2f\ud v\right)=0}
\end{cases}
\end{equation}
\end{prop}
Finally, under the same assumptions on $B$, Boltzmann's H-theorem is true:

\begin{prop}
 Let the collision kernel $B$ satisfy \eqref{hyp1B} and \eqref{hypsurB}. Assume that the initial datum satisfies \eqref{hypf00}, and
\[
 \iint_{\Omega\times\R^3}(|x|^2+|v|^2)f_0(x,v)\ud x\ud v <+\infty,
\]
and let $f$ be the solution to \eqref{BBFD}-\eqref{boundary}. Then, 
\[
 \iint_{\Omega\times\R^3}s(f)(t)\ud v\ud x +\int_0^t \int_\Omega D(f)\ud x\ud s=\iint_{\Omega\times\R^3}s(f_0).
\]
\end{prop}

\subsection{The linearized collision operator and the fluxes}

It is sometimes interesting to linearize the Fermi-Dirac collision operator $Q(f)$ around a global equilibrium. Let $P_\e=\frac{M_\e}{1+M_\e}$ be a global equilibrium depending on the parameter $\e$, with $M_\e = a_\e
e^{-\frac{|v-u_\e|^2}{2b_\e}}$. Then the linearized collision operator writes

\bee
\begin{split}
 \mc L_\eps g&=\int_{S^2}\int_{\R^3}B(v-v_*,\omega)\frac{M_{\eps,*}}{1+M_{\eps,*}}\frac{1}{1+M_\eps'}\frac{1}{1+M_{\eps,*}'}\\
&\quad\qquad\times\left(g(1+M_\eps)+g_*(1+M_{\eps,*})-g'(1+M_\eps')-g_*'(1+M_{\eps,*}')\right)\ud v_*\ud\omega.
\end{split}
\eee

The properties of this operator will be detailed and proved in Section \ref{intermediate}. Let us just mention that if the collision kernel $B$ satisfies 
Grad's cutoff assumptions \eqref{Grad1}-\eqref{Grad2} then $\mc L_\e$ is a Hilbert-Schmidt operator with kernel
 \[
 \textrm{Ker} \mc L_\eps = \textrm{span}\left\{\frac{1}{1+M_\eps},\frac{v}{1+M_\eps},\frac{|v|^2}{1+M_\eps}\right\}.
\]
 This property will be crucial in the proof of our main result.

\begin{rem}
Let us make here an important remark. The collision kernel $B$ cannot satisfy both assumptions \eqref{hyp1B} (required to prove the existence of solutions) and \eqref{Grad1}-\eqref{Grad2}. 
We will work with two different collision kernels : $B_\e(z,\omega) =\mathbbm{1}_{|z|\leq \frac{1}{\eps^2}} $ satisfies \eqref{hyp1B} and will be used in the definitions of the collision operator $Q(f_\e)$ and the entropy dissipation $D(f_\e)$.

However, $B= \mathbbm{1}$ will be considered when studying the properties of the linear operator.

In proving the hydrodynamic limit, the collision kernel $B_\e$ will naturally appear in expressions to control. When the linear operator $\mc L_\e$ will be needed, we will write

\[
B_\e = B + (B_\e -B)
\]
and use the linear operator properties on the term containing $B$, while the other one will be controlled by hands using the fact that $B_\e$ and $B$ are equal on a very large ball.

\end{rem}

Let us define the energy and heat fluxes
\[
 \Phi_\eps=\frac{1}{b_\eps}\left((v-\eps u)^{\otimes 2}-|v-\eps u|^2 I\right)
\]
\[
 \Psi_\eps=\frac{v-\eps u}{b_\eps^2}\left(|v-\eps u|^2-\tau_\eps\right),
\]
where $\tau_\eps$ is the unique number such that $\frac{\Psi_\eps}{1+M_\eps}\in (\textrm{Ker}\mc L_\eps)^\bot$ and $I$ is the $2\times 2$ identity matrix. These quantities naturally appear in Section \ref{sec:modulated}. Taking $\e=0$, we can prove that
\[
\Phi, \Psi \in (\textrm{Ker}\mc L_0)^\bot
\]
where $\Phi=\Phi_0$ and $\Psi=\Psi_0$. Since $\mc L_0$ is a Hilbert-Schmidt operator, we can define

\be
\label{poly}
 \tilde \Phi = \mc L_0^{-1}\Phi\quad \textrm{and}\quad \tilde \Psi = \mc L_0^{-1}\Psi.
\ee

$\tilde\Phi$ and $\tilde \Psi$ will appear in the control of the flux terms.

\section{Strategy and main result}
\label{sec:strategy}
%In this Section, we explain the method of proof, give our main theorem and some intermadiate results.

%In order to show that with the scaling under consideration, the solutions $f_\eps$ to equation \eqref{scaled} will converge to some
%equilibrium solution the parameters of which satisfy the incompressible Euler equations, we use the relative entropy method 

In investigating the incompressible Euler limit of the Boltzmann-Fermi equation, we use the modulated entropy method as in the
classical case \cite{saintray}. Let us first introduce the relative entropy of the solution $f_\e$ to \eqref{scaled} with respect to the
global thermodynamic equilibrium $P_0$ (defined by \eqref{P_0}):
\[
\begin{split}
  H(f_\eps|P_0)=&\int_{\Omega\times\R^3}\bigg(f_\eps \log \frac{f_\eps}{P_0}-f_\eps+P_0\\
&\qquad\qquad+(1-f_\eps)\log\frac{1-f_\eps}{1-P_0}-(1-f_\eps)+1-P_0\bigg)\ud v\ud x.
\end{split}
\]
It is a non-negative Lyapunov functional for the Boltzmann-Fermi equation thanks to the H theorem. In the fast relaxation limit, $f_\e$ is
supposed to be close to a thermodynamic equilibrium, we therefore define the modulated entropy by

\[
\begin{split}
  H(f_\eps|P_\eps)=&\int_{\Omega\times\R^3}\bigg(f_\eps \log \frac{f_\eps}{P_\eps}-f_\eps+P_\eps\\ &\qquad\qquad+(1-f_\eps)\log\frac{1-f_\eps}{1-P_\eps}-(1-f_\eps)+1-P_\eps\bigg)\ud v\ud x.
\end{split}
\]
where $P_\eps =\frac{M_\eps}{1+M_\eps}$ is some local equilibrium which approximates the solution $f_\eps$, with $M_\e = a_\e
e^{-\frac{|v-u_\e|^2}{2b_\e}}$. The relative entropy functional
measures this approximation since
\[
 H(f_\eps|P_\eps)\geq \int_{\Omega\times\R^3}\left(\sqrt{f_\eps}-\sqrt{P_\eps}\right)^2\ud v\ud x,
\]
consequence of the pointwise inequality
\[
 (\sqrt{x}-\sqrt{y})^2\le x\log\frac{x}{y}-x+y\qquad \forall x,y>0.
\]
Our strategy consists in studying the time evolution of $\frac{1}{\e^2}H(f_\e|P_\e)$. This has the good scaling since we want to observe the fluctuations of $f_\e$ around an equilibrium. To do that we 
\begin{itemize}
  \item determine the equation satisfied by $\frac{1}{\e^2}H(f_\e|P_\e)$, using the local conservation laws satisfied by $f_\e$; it contains an acoustic (fast oscillating) term and flux terms
  \item control the acoustic terms by specifying $a_\e$, $u_\e$ and $b_\e$
  \item control the flux terms thanks to a priori bounds on the linear collision operator.
\end{itemize}

Our main result is the following:
\begin{thm}
\label{limite}
Let $\Omega$ be some regular bounded domain of $\R^3$. Let $(f_{\eps,in})$ be a family of measurable nonnegative functions over
$\Omega\times \R^3$ satisfying the bounds
\[
 \int_{\Omega}\int_{\R^3}(1+|x|^2+|v|^3)f_{\eps,in}\ud v\ud x \leq C_\eps,
\]
\[
 0\leq f_{\eps,in}\leq 1,
\]
and the scaling condition 
\be
\label{scal}
 \frac{1}{\eps^2}H(f_{\eps,in}|P_{0})\leq C.
\ee
Without loss of generality, assume that the fluctuation $g_{\eps,in}$ defined by $f_{\eps,in}=P_0(1+\eps g_{\eps,in})$ converges
\[
 g_{\eps,in} \to g_{in}=\frac{1}{1+M_0}\left(\frac{a_{in}}{a_0} +\frac{u_{in}}{b_0} \cdot v + \frac{b_{in}}{2b_0}|v|^2\right).
\]
We assume moreover that the initial data is well-prepared: 
\be
\label{bienprepare}
 \frac{1}{\eps^2}H(f_{\eps,in}|P_{\eps,in})\to 0.
\ee
Let $f_\eps$ be some family of solutions to the scaled Boltzmann-Fermi equation

\be
\begin{cases}
&\eps\partial_t f_\eps + v.\nabla_x f_\eps = \frac{1}{\eps^q}Q(f_\eps)\\
&f_\eps(0,x,v) = f_{\eps,in}(x,v) \quad \textrm{on}~\Omega\times\R^3 
\end{cases}
\ee
with $q>1$, and with a maxwellian truncated collision kernel $B_\eps(z,\omega)=\mathbbm{1}_{|z|\leq \frac{1}{\eps^2}}$, endowed with the boundary condition
\be
\label{boundary2}
 f_\eps(t,x,v)=f_\eps(t,x,R_x(v))\qquad\forall (x,v)\in \partial\Omega\times\R^3 ~~\textrm{such that}~~n(x).v<0
\ee
where $R_x(v)$ is the specular reflection law \eqref{specular}. Assume that we have the following control on the tails :
\be
\label{control}
\int_{\R^3}P_0\left(\frac{f_\eps-P_0}{P_0}\right)^2\ud v \leq C ~~\textrm{a.e.}~
\ee
and, finally, assume that 
\be
\label{assum:poly}
 \tilde \Phi = \mc L_0^{-1}\Phi\quad \textrm{and}\quad \tilde \Psi = \mc L_0^{-1}\Psi\quad \textrm{are at most polynomial as}\quad|v|\to \infty.
\ee
Then the fluctuation $(g_\e)$ defined by $f_\e=P_0(1+\e g_\e)$ converges in $L^1([0,T];(1+|v|^2)P_0\ud v\ud x)$ weak 
\[
 g_\e \underset{\e\to 0}{\longrightarrow}\frac{1}{1+M_0}\left( \frac{\bar a}{a_0}+\frac{\bar u}{b_0} \cdot v +\frac{\bar b}{2b_0}|v|^2\right)
\]
where $(\bar a, \bar u, \bar b)$ is the unique Lipschitz solution to the incompressible Euler equations
\[
 \partial_t \bar u + \bar u.\nabla_x \bar u + \nabla_x p=0,\qquad\nabla_x.\bar u=0
\]
supplemented with
\[
 \partial_t \bar a + \bar u.\nabla_x \bar a=0,\qquad\qquad\partial_t \bar b+\bar u.\nabla_x \bar b=0,\qquad\qquad\nabla_x(b_0\bar a+\tau_0\bar b)=0
\]
on $[0,T]$, and $T$ is the maximal lifespan of the solution.
\end{thm}
%
%
%% In \eqref{assum:poly}, $\mc L_0$ is the linearized collision operator around $P_0$. We will see in the next subsection that it is a
%% Hilbert-Schmidt operator and that $\Phi$ and $\Psi$ are orthogonal to its kernel. This gives a sense to $\tilde \Phi$ and $\tilde \Psi$. 

The result is not very surprising since it is very similar to the classical case \cite{saintray}. Assumption \eqref{assum:poly} was proved in the classical case \cite{polynomial}, and it should be possible (although very technical) to prove it also in the fermion case.

Assumption \eqref{bienprepare} allows us to consider solutions which are almost instantaneously at thermodynamic equilibrium, thus avoiding the Knudsen layer. To get rid of it we need to get a better understanding of the relaxation mecanism.

As to assumption \eqref{control}, it is a purely technical control of large velocities which allows us to treat flux terms. The exact same assumption is done in the classical case and it seems to be a great challenge to avoid it. It has been proved that some regular solutions to the classical Boltzmann equation satisfying \eqref{control} exist, but such solutions in the fermion case are still missing to the author's knowledge.

The cubic nonlinearity of the collision operator adds a lot of technical difficulties in the proof: a lot more terms, respect to the classical case, appear and must be controlled. The $L^\infty$ bound 
is of great help to treat these difficulties, but does not bring more, so that at the end we do not recover a better result than in the classical case.

%Our strategy of proof is basically the same than for the classical case \cite{saintray}. The main difference is in the control of the flux terms. %We still use a decomposition of the linearized collision operator $\mc L_\eps$, but it contains more terms.

\section{The modulated entropy}
\label{sec:modulated}

The main idea to prove the convergence of the solutions to the scaled Boltzmann-Fermi equation toward a solution of the incompressible Euler equation is to study the time evolution of the modulated entropy. Let 
\[
     P_\eps=\frac{M_\eps}{1+M_\eps},\qquad  M_\eps=a_\eps e^{-\frac{|v-\eps u|^2}{2b_\eps}},
\]
be a local equilibrium which approximates the solution $f_\eps$, and take
\be
\label{ab}
 a_\eps=  e^{a_0+\eps a_1},\qquad\qquad b_\eps=b_0e^{\eps b_1}.
\ee
Recall that $a_0$ and $b_0$ are positive constants which were defined in Subsection \ref{relation}; $a_1$, $u$ and $b_1$ are regular functions of $(t,x)$. 
%% The modulated entropy is defined by
%% \[
%%  H(f_\eps|P_\eps)=\int_{\Omega\times\R^3}\left(f_\eps \log \frac{f_\eps}{P_\eps}-f_\eps+P_\eps +(1-f_\eps)\log\frac{1-f_\eps}{1-P_\eps}-(1-f_\eps)+1-P_\eps\right)\ud v\ud x.
%% \]

\begin{prop}
\label{modulated entropy identity}
 Any solution to the scaled Boltzmann-Fermi equation \eqref{scaled} satisfies the following identity
\bee
\begin{split}
\frac{1}{\eps^2}H&(f_\eps|P_\eps)(t)=\frac{1}{\eps^2}H(f_{\eps,in}|P_{\eps,in}) -\frac{1}{\eps^{q+3}}\int_0^t\int_{\Omega}D(f_\eps)\ud x\ud s\\
&-\frac{1}{\eps}\int_0^t\iint_{\Omega\times\R^3}f_\eps\left(1,\frac{v-\eps u}{b_\eps}, \frac{|v-\eps u|^2}{2b_\eps}\right).\mc A_\eps(a_1,u,b_1)\ud x\ud v \ud s\\
&-\frac{1}{2\eps^2}\int_0^t\iint_{\Omega\times\R^3}f_\eps b_\eps\Psi_\eps.\nabla_xb_1\ud x\ud v\ud s -\frac{1}{2\eps^2}\int_0^t\iint_{\Omega\times\R^3}f_\eps (\Phi_\eps\cdot\nabla_x u)\ud x\ud v\ud s\\
&+\frac{1}{3\eps^2}\int_0^t\left[\frac{\ud}{\ud t}\iint_{\Omega\times\R^3}\frac{e^{-\eps b_1}}{b_0}|v|^2\frac{e^{a_0+\eps a_1-\frac{e^{-\eps b_1}|v|^2}{2b_0}}}{1+e^{a_0+\eps a_1-\frac{e^{-\eps b_1}|v|^2}{2b_0}}}\ud v\ud x \right]\ud s.
\end{split}
\eee
where the acceleration operator $\mc A_\eps(a_1,u,b_1)$ is defined by 
\be
\label{acceleration}
 \mc A_\eps(a,u,b)=\left(\begin{array}{l}
                    \partial_t a+u.\nabla_x a \\
		    \partial_t u +u.\nabla_x u+ b_0\frac{e^{\eps b}-1}{\eps} \nabla_x a +\frac{1}{\eps}(\tau_\eps-\tau_0)\nabla_x b+\frac{1}{\eps} \nabla_x \left( b_0a +\tau_0 b\right)\\
		    \partial_t b +u.\nabla_x b+\frac{2}{3\eps}\nabla_x. u
                   \end{array}\right)
\ee
and $\tau_0$ is the limit of $\tau_\e$ as $\e$ goes to $0$. 
\end{prop}

\begin{proof}

The modulated entropy can be written as
\[
 H(f_\eps|P_\eps)=H(f_\eps|P_0)+\int_{\Omega\times\R^3}\left(f_\eps\log\frac{M_0}{M_\eps}+\log\frac{1+M_\eps}{1+M_0}\right)\ud v\ud x.
\]
The H-theorem allows to write
\be
\label{entropymod}
\begin{split}
H(f_\eps|P_\eps)(t)&=H(f_{\eps,in}|P_{\eps,in}) -\frac{1}{\eps^{q+1}}\int_0^t\int_{\Omega}D(f_\eps)\ud x\ud s\\
&+\int_0^t\left[\frac{\ud}{\ud t}\iint_{\Omega\times\R^3}f_\eps\left(-\log a_\eps+\frac{|v-\eps u|^2}{2b_\eps}-\frac{|v|^2}{2}\right)\ud x\ud v \right]\ud s\\
&+\int_0^t\left[\frac{\ud}{\ud t}\iint_{\Omega\times\R^3}\log\left(1+M_\eps\right)\ud x\ud v \right]\ud s.
\end{split}
\ee
An integration per part changes the last term:
\[
 \iint_{\Omega\times\R^3}\log\left(1+M_\eps\right)\ud x\ud v=\frac{1}{3}\iint_{\Omega\times\R^3}\frac{1}{b_\eps}|v|^2\frac{ a_\e e^{-\frac{|v|^2}{2b_\eps}}}{1+a_\e e^{-\frac{|v|^2}{2b_\eps}}}\ud v\ud x.
\]
Computing the time derivative and using the conservation laws \eqref{conslaws} leads to:

%% \bee
%% \begin{split}
%% &\frac{\ud}{\ud t}\iint_{\Omega\times\R^3}f_\eps\left(-\log a_\e+\frac{|v-\eps u|^2}{2b_\eps}-\frac{|v|^2}{2}\right)\ud x\ud v\\
%% &=\iint_{\Omega\times\R^3}\bigg[\partial_tf_\eps\left(-\log a_\e+\frac{|v-\eps u|^2}{2b_\eps}-\frac{|v|^2}{2}\right)\\
%% &+f_\eps\left(-\partial_t \log a_\e -\frac{\partial_t b_\eps}{b_\eps}\frac{|v-\eps u|^2}{2b_\eps}-\frac{\eps}{b_\eps}(v-\eps u) \cdot \partial_t u\right)\bigg]\ud x\ud v.
%% \end{split}
%% \eee
%% Using the conservation laws \eqref{conslaws}, it comes

\bee
\begin{split}
 \frac{\ud}{\ud t}\iint_{\Omega\times\R^3}&f_\eps\left(-\log a_\e+\frac{|v-\eps u|^2}{2b_\eps}-\frac{|v|^2}{2}\right)\ud x\ud v\\
&=-\iint_{\Omega\times\R^3}f_\eps\bigg[\partial_t \log a_\e+u \cdot \nabla_x \log a_\e \\
&~~~~~+\frac{\eps}{b_\eps}(v-\eps u)\left(\partial_t u +u \cdot \nabla_x u+\frac{1}{\eps}\nabla_x u:(v-\eps u)+\frac{1}{\eps^2}b_\eps\nabla_x \log a_\e\right)\\
&~~~~~+\frac{|v-\eps u|^2}{2b_\eps}\left(\frac{\partial_t b_\eps}{b_\eps} +u\cdot \frac{\nabla_x b_\eps}{b_\eps}+\frac{1}{\eps}(v-\eps u) \cdot \frac{\nabla_x b_\eps}{b_\eps} \right)\bigg]\ud v\ud x.
\end{split}
\eee
Introducing the rescaled translated versions of the momentum and energy fluxes

\[
 \Phi_\eps=\frac{1}{b_\eps}\left((v-\eps u)^{\otimes 2}-|v-\eps u|^2 I\right)~~~~\textrm{and}~~~~ \Psi_\eps=\frac{v-\eps u}{b_\eps^2}\left(|v-\eps u|^2-\tau_\eps\right)
\]
and replacing $a_\e$ and $b_\e$ by \eqref{ab} we obtain:

%% \bee
%% \begin{split}
%%  \frac{\ud}{\ud t}\iint_{\Omega\times\R^3}&f_\eps\left(-\log a_\e+\frac{|v-\eps u|^2}{2b_\eps}-\frac{|v|^2}{2}\right)\ud x\ud v\\
%% &=-\iint_{\Omega\times\R^3}f_\eps\bigg[\partial_t \log a_\e+u\cdot\nabla_x \log a_\e \\
%% &~~~~~+\frac{\eps}{b_\eps}(v-\eps u)\left(\partial_t u +u\cdot\nabla_x u+\frac{1}{\eps^2}b_\eps \nabla_x \log a_\e +\frac{1}{\eps^2}\tau_\eps\frac{\nabla_x b_\eps}{b_\eps}\right)\\
%% &~~~~~+\frac{|v-\eps u|^2}{2b_\eps}\left(\frac{\partial_t b_\eps}{b_\eps} +u\cdot\frac{\nabla_x b_\eps}{b_\eps} +\frac{2}{3\e}\nabla_x\cdot u\right)\bigg]\ud v\ud x\\
%% &~~~~-\frac{1}{2\eps}\iint_{\Omega\times\R^3}f_\eps\Psi_\eps\cdot\nabla_xb_\eps\ud x\ud v -\frac{1}{2\e}\iint_{\Omega\times\R^3}f_\eps (\Phi_\eps\cdot\nabla_x u)\ud x\ud v.
%% \end{split}
%% \eee
%% We now replace $a_\e$ and $b_\e$ by \eqref{ab}, which leads to

\bee
\begin{split}
 \frac{\ud}{\ud t}\iint_{\Omega\times\R^3}&f_\eps\left(-\log a_\e+\frac{|v-\eps u|^2}{2b_\eps}-\frac{|v|^2}{2}\right)\ud x\ud v\\
&=-\iint_{\Omega\times\R^3}f_\eps\bigg[\eps  \left(\partial_t a_1+u\cdot\nabla_x a_1 \right)\\
&~~~~~+\frac{\eps}{b_\eps}(v-\eps u)\left(\partial_t u +u\cdot\nabla_x u+\frac{1}{\eps} b_\eps \nabla_x a_1 +\frac{1}{\eps}\tau_\eps\nabla_x b_1\right)\\
&~~~~~+\eps\frac{|v-\eps u|^2}{2b_\eps}\left(\partial_t b_1 +u\cdot\nabla_x b_1+\frac{2}{3\eps}\nabla_x\cdot u \right)\bigg]\ud v\ud x\\
&~~~~-\frac{1}{2\eps}\iint_{\Omega\times\R^3}f_\eps\Psi_\eps\cdot\nabla_xb_\eps\ud x\ud v -\frac{1}{2\e}\iint_{\Omega\times\R^3}f_\eps (\Phi_\eps\cdot\nabla_x u)\ud x\ud v.
\end{split}
\eee
Summarizing, the modulated entropy satisfies 
\bee
\begin{split}
\frac{1}{\eps^2}H(f_\eps|P_\eps)(t)&=\frac{1}{\eps^2}H(f_{\eps,in}|P_{\eps,in}) -\frac{1}{\eps^{q+3}}\int_0^t\int_{\Omega}D(f_\eps)\ud x\ud s\\
&-\frac{1}{\eps}\int_0^t\iint_{\Omega\times\R^3}f_\eps\left(1,\frac{v-\eps u}{b_\eps}, \frac{|v-\eps u|^2}{2b_\eps}\right)\cdot\mc A_\eps(a_1,u,b_1)\ud x\ud v \ud s\\
&-\frac{1}{2\eps^2}\int_0^t\iint_{\Omega\times\R^3}f_\eps\left( b_\eps\Psi_\eps\cdot\nabla_xb_1+ \Phi_\eps\cdot\nabla_x u\right)\ud x\ud v\ud s\\
&+\frac{1}{3\eps^2}\int_0^t\left[\frac{\ud}{\ud t}\iint_{\Omega\times\R^3}\frac{e^{-\eps b_1}}{b_0}|v|^2\frac{e^{a_0+\eps a_1-\frac{e^{-\eps b_1}|v|^2}{2b_0}}}{1+e^{a_0+\eps a_1-\frac{e^{-\eps b_1}|v|^2}{2b_0}}}\ud v\ud x \right]\ud s.
\end{split}
\eee
where the acceleration operator is defined by \eqref{acceleration}.

\end{proof}

\section{Construction of the approximate solutions}
\label{sec:construction}

A global equilibrium solution is not expected to be a good approximation of $f_\e$ in the fast relaxation limit since fast oscillations can take place, such as acoustic waves. Hence we have to find correctors in order to obtain a refined approximation which will lead to the convenient asymptotics.

We want to find $V_\eps=(a_1^\eps,u_\eps,b_1^\eps)$ such that $\mc A_\eps(a_1^\eps,u_\eps,b_1^\eps)\to 0$ in $L^2$, that is, such that $V_\e$ is an approximate solution of the system
\be
\label{eqV}
 \partial_t V + \frac{1}{\eps}WV+B(V,V)= 0,
\ee
where 
\[
 V=(a_1,u,b_1),
\]
%% the penalization operator
\[
 WV=\Big(0,\nabla_x \left( b_0a_1 +\tau_0 b_1\right),\frac{2}{3}\nabla_x. u\Big)
\]
%% and the bilinear operator
and
\[
 B(V,V)=\left(\begin{array}{l}
                u.\nabla_x a_1\\
		u.\nabla_x u+ b_0b_1 \nabla_x a_1 +\tau_1\nabla_x b_1\\
		u.\nabla_x b_1
              \end{array}
		\right)
\]
with moreover the constraint
\be
\label{constraint}
 A_\eps=\frac{1}{\eps^2}\frac{\ud}{\ud t}\iint_{\Omega\times\R^3}\frac{e^{-\eps b_1}}{b_0}|v|^2\frac{e^{a_0+\eps a_1-\frac{e^{-\eps b_1}|v|^2}{2b_0}}}{1+e^{a_0+\eps a_1-\frac{e^{-\eps b_1}|v|^2}{2b_0}}}\ud v\ud x\to 0.
\ee

\begin{rem}
B is a bilinear operator since we can decompose $\tau_\e$ as $\tau_\e = \tau_0 + \e \tau_1 + o(\e)$ and $\tau_1$ is a linear combination of $a_1$ and $b_1$.
\end{rem}

Such a constrution is done by filtering methods, and we refer to \cite{saintray} for all the details. We will just here mention the main points of the proof of

\begin{thm}
 \label{approxsols}
Let $(a_1^{in}, u^{in}, b_1^{in})$ belong to $H^s(\Omega)$ for some $s>\frac{5}{2}$. Then there exists some $T>0$ and some family $(a_1^{\eps,N}, u^{\eps,N}, b_1^{\eps,N})$ such that 
\[
 \sup_{N\in\N} \lim_{\eps\to 0} \|(a_1^{\eps,N}, u^{\eps,N}, b_1^{\eps,N})\|_{L^\infty([0,T];H^s(\Omega))}\leq C_T,
\]

\[
 (a_1^{\eps,N,in}, u^{\eps,N,in}, b_1^{\eps,N,in})\to (a_1^{in}, u^{in}, b_1^{in})\quad \textrm{in} ~H^s(\Omega)~\textrm{as}~\eps\to 0 ~\textrm{and}~N\to +\infty,
\]

\[
\mc A_\eps  (a_1^{\eps,N}, u^{\eps,N}, b_1^{\eps,N}) \to 0 \quad \textrm{in} ~L^2([0,T]\times\Omega)~\textrm{as}~\eps\to 0 ~\textrm{and}~N\to +\infty,
\]
and
\[
  A_\eps^N=\frac{1}{\eps^2}\frac{\ud}{\ud t}\iint_{\Omega\times\R^3}\frac{e^{-\eps b_1^{\eps,N}}}{b_0}|v|^2\frac{e^{a_0+\eps a_1^{\eps,N}-\frac{e^{-\eps b_1^{\eps,N}}|v|^2}{2b_0}}}{1+e^{a_0+\eps a_1^{\eps,N}-\frac{e^{-\eps b_1^{\eps,N}}|v|^2}{2b_0}}}\ud v\ud x\longrightarrow 0
\]
in $L^1([0,T])$ as $\eps\to 0$ and $N\to +\infty$.

\end{thm}

%% for which is the mass conservation is required, as shown by Proposition \ref{A tend vers 0}.

In $L^2$ equipped with the norm
\[
 \|V\|^2_{L^2}=\|a_1\|^2_{L^2}+\|u\|^2_{L^2}+\frac{3}{2\tau_0}\|b_0 a_1 + \tau_0 b_1\|^2_{L^2},
\]
the operator $W$ is skew-symmetric. We can thus define the associated semigroup $\mc W$. If we conjugate \eqref{eqV} with $\mc W(\frac{t}{\eps})$, it comes
\[
 \partial_t (\mc W \left(\frac{t}{\eps}\right) V) + \mc W \left(\frac{t}{\eps}\right)B(V,V)=0.
\]
Defining
\[
 \tilde V= \mc W\left(\frac{t}{\eps}\right) V,
\]
it becomes
\be
\label{acoustic}
 \partial_t \tilde V +\mc W \left(\frac{t}{\eps}\right)B\left(\mc W \left(-\frac{t}{\eps}\right)\tilde V,\mc W \left(-\frac{t}{\eps}\right)\tilde V\right)=0.
\ee
We obtain the first order approximation by taking strong limits in the filtered systems. However, the error is not expected to converge strongly to 0 due to high frequency oscillations. We therefore need to construct a second order approximation. We then need a third order approximation in order to have mass conservation at sufficient order.

\subsection{Study of $W$}
We expect the solutions of \eqref{acoustic} to have a very different behaviour depending on the spectrum of $W$. Since it is skew-symmetric, its eigenvalues are purely imaginary, and its eigenvectors satisfy
\[
 WV_\lambda=i\lambda V_\lambda,
\]
which implies if $\lambda \neq 0$
%% \[
%%  \begin{cases}
%%   &0=i\lambda a_\lambda\\
%%   &\nabla_x (b_0 a_\lambda +\tau_0 b_\lambda)=i\lambda u_\lambda\\
%%   &\frac{2}{3}\nabla_x.u_\lambda=i\lambda b_\lambda
%%  \end{cases}
%% \]
%% In other words, if $\lambda \neq 0$, the eigenvectors satisfy
\[
 \begin{cases}
  &a_\lambda=0\\
  &\Delta_x b_\lambda=-\frac{3}{2\tau_0}\lambda^2 b_\lambda\\
  & \nabla_x b_\lambda=i\frac{\lambda}{\tau_0}u_\lambda
 \end{cases}.
\]
Hence, the operator $W$ has the same spectral structure as the laplacian on $\Omega$, which means that is is diagonalizable on the orthogonal of its kernel.%%  Indeed, there exists an hilbertian basis $(\phi_n)_{n\in I}$ of $L^2(\Omega)$ such that
%% \[
%%  -\Delta_x \phi_n =\mu_n \phi_n \quad \textrm{with}\quad n.\nabla_x\phi_n=0\quad \textrm{on}\quad \partial\Omega
%% \]
%% with $\mu_n >0$ for $n\in \N^*$ and
%% \[
%%  \Delta_x \phi_n=0\quad \textrm{with}\quad n.\nabla_x\phi_n=0\quad \textrm{on}\quad \partial\Omega
%% \]
%% for $n\in I\setminus \N^*$. This yields that the family $(\psi_n)_{n\in \Z}$ defined by
%% \[
%%  \psi_n=\left(0,i ~\textrm{sgn}(n)\sqrt{\frac{3\tau_0}{2\mu_{|n|}}}\nabla_x \phi_{|n|}, \phi_{|n|}\right)
%% \]
%% is an hilbertian basis of $(\textrm{ker }W)^\bot$.
An explicit computation shows that the orthogonal projection on $\textrm{ker }W$ is
\[
 \Pi_0 (a,u,b)=\left(a, Pu, \frac{|\Omega|^{-1}\int(b_0 a + \tau_0 b)\ud x-b_0 a }{\tau_0}\right)
\]
where $P$ is the Leray projection, that is, the orthogonal projection onto the divergence-free vectors. In the sequel, $\Pi_\lambda$ is the orthogonal projection onto $\textrm{ker}(W-\lambda I)$.

\subsection{Construction}
We just sketch here the construction of approximate solutions to \eqref{acoustic}, and we refer to \cite{saintray} for all the details. The method is as follows:

\begin{itemize}
 \item Decompose the operator $W$ as $W=\sum_{\lambda\in  \sigma_p}\lambda \Pi_\lambda$ and plug this into \eqref{acoustic}. It leads to
\[
 \partial_t \tilde V + \sum_{k_1, k_2, k_3 \in \sigma_p}e^{i\frac{t}{\eps}(\lambda_{k_1}-\lambda_{k_2}-\lambda_{k_3})}\Pi_{\lambda_{{k_1}}}B(\Pi_{\lambda_{k_2}}\tilde V,\Pi_{\lambda_{k_3}}\tilde V)=0
\]
 \item The first order approximation $\tilde V_0$ of $\tilde V$ is defined as the solution of equation \eqref{acoustic} when we take into account only the resonant modes, that is,
\[
  \partial_t \tilde V_0 + \sum_{\lambda_{k_1}= \lambda_{k_2}+ \lambda_{k_3} }\Pi_{\lambda_{k_1}}B(\Pi_{\lambda_{k_2}}\tilde V_0,\Pi_{\lambda_{k_3}}\tilde V_0)=0
\]
This equation is known to have solutions in $L^\infty_{loc}([0,T_*),H^s(\Omega))$ provided that $V^{in}\in H^s(\Omega)$ and $s>\frac{5}{2}$. However, $\tilde V_0$ is not an approximation of $\tilde V$ in the sense that $\mc A_\eps  (\tilde V_0)$ converges weakly but not strongly in $L^2$ to 0. We therefore have to add correctors.

 \item Hence, we construct the second order approximation plugging $\tilde V=\tilde V_0 +\eps \tilde V_1+o(\eps)$ into \eqref{acoustic}. To avoid the problem of small divisors, we introduce the projection onto a finite dimensional subset of $C^\infty (\bar \Omega)$ $J_N$, which is the orthogonal projection onto the $N$ first harmonics of $W$ and the $N$ first harmonics in $\textrm{ker } W$. We then denote $\tilde V_0^N=J_N\tilde V_0$, and we check that $\int \tilde V_0 \ud x = \int \tilde V_0^N\ud x$. We define $\tilde V_1^N$ by
\[
 \tilde V_1^N=J_N \sum_{\lambda_{k_1}+\lambda_{k_2}\neq \lambda_{k}} \frac{\exp\left(\frac{it}{\eps}(\lambda_k-\lambda_{k_1}-\lambda_{k_2})\right)} {i(\lambda_{k_1}+\lambda_{k_2}-\lambda_k)}\Pi_{\lambda_k}B(\Pi_{\lambda_{k_1}}\tilde V_0^N,\Pi_{\lambda_{k_2}}\tilde V_0^N).
\]
Then, $\tilde V_0^N+\eps \tilde V_1^N$ is an approximate solution of \eqref{acoustic} strongly in $L^2$. However, it does not satisfy $\frac{1}{\eps^2}\int (\tilde V_0^N+\eps \tilde V_1^N)\ud x \to 0$ when $\eps\to 0$. 

\item We therefore define the third order approximation $\tilde V_2^N$ as in \cite{saintray} by plugging $\tilde V_\eps^N=\tilde V_0^N+\eps \tilde V_1^N +\eps^2 \tilde V_2^N+o(\e^2)$ into \eqref{acoustic}, and we check that
	\begin{itemize}
	\item $\tilde V_\eps^N=\tilde V_0^N+\eps \tilde V_1^N +\eps^2 \tilde V_2^N$ is an approximate solution of \eqref{acoustic},
	\item  $\frac{1}{\eps^2}\int \tilde V_\eps^N\ud x \to 0$ when $\eps\to 0$.
	\end{itemize}

 %\item We get an approximate solution $\tilde V^\eps$ of equation \eqref{acoustic} such that $\tilde V^\eps \to \tilde V_0$ as $\eps$ go to 0 and such that $\partial_t \int \tilde V^\eps \ud x =0$.
\end{itemize}

The key ingredient for the above computations to work (see \cite{saintray}) is that, for $\lambda,\mu \neq 0$ with $\lambda \neq \mu$, we have
\[
 \Pi_0 B(\Pi_\lambda V,\Pi_\mu V)=0.
\]
Indeed, writing $\Pi_\lambda V=(a_\lambda, u_\lambda, b_\lambda)$, we get
\bee
\begin{split}
  \Pi_0 B(\Pi_\lambda V,\Pi_\mu V)=\frac{1}{2}\Pi_0 \begin{pmatrix}
                                                     u_\lambda.\nabla_x a_\mu + u_\mu.\nabla_x a_\lambda \\
						     u_\lambda.\nabla_x u_\mu + u_\mu.\nabla_x u_\lambda +b_0 b_\lambda\nabla_x a_\mu + b_0 b_\mu\nabla_x a_\lambda \\
						     \qquad\qquad+ (c^1a_\lambda + c^2b_\lambda)\nabla_x b_\mu +(c^1 a_\mu + c^2b_\mu)\nabla_x b_\lambda\\
						     u_\lambda.\nabla_x b_\mu+u_\mu.\nabla_xb_\lambda
                                                    \end{pmatrix};
\end{split}
\eee
but $a_\lambda=a_\mu=0$ since $\lambda,\mu\neq0$, so that
\bee
\begin{split}
  \Pi_0 B(\Pi_\lambda V,\Pi_\mu V)&=\frac{1}{2}\Pi_0 \begin{pmatrix}
                                                     0 \\
						     u_\lambda.\nabla_x u_\mu + u_\mu.\nabla_x u_\lambda +  c^2b_\lambda\nabla_x b_\mu+c^2b_\mu\nabla_x b_\lambda\\
						     u_\lambda.\nabla_x b_\mu+u_\mu.\nabla_xb_\lambda
                                                    \end{pmatrix}\\
&=\begin{pmatrix}
   0\\
   P(u_\lambda.\nabla_x u_\mu + u_\mu.\nabla_x u_\lambda)\\
   \oint(u_\lambda.\nabla_x b_\mu + u_\mu . \nabla_x b_\lambda)\ud x)
  \end{pmatrix}
.
\end{split}
\eee
Note that $\int u_\lambda.\nabla_x b_\mu \ud x = -\frac{3\mu^2}{2i\lambda}\int b_\lambda.b_\mu\ud x=0$ since $\lambda\neq\mu$. Moreover,
\[
\left( u_\lambda.\nabla_x u_\mu + u_\mu.\nabla_x u_\lambda\right)_i=\left(\nabla_x(u_\lambda.u_\mu)\right) +u_\lambda^j\left(\partial_ju_\mu^i -\partial_i u_\mu^j\right)+u_\mu^j\left(\partial_j u_\lambda^i-\partial_i u_\lambda^j\right),
\]
so that
\[
 \Pi_0 B(\Pi_\lambda V,\Pi_\mu V)=0.
\]
In addition, if $\lambda\neq 0$, 
\[
 \int_{\Omega} b_\lambda\ud x =-i\frac{2}{3\lambda}\int_{\Omega} \nabla_x.u_\lambda \ud x = i\frac{2}{3\lambda}\int_{\partial \Omega} u.n \ud \sigma_x=0.
\]
At the end, we get that
\[
 \frac{\ud}{\ud t} \int_{\Omega}\tilde b\ud x=\frac{\ud}{\ud t} \int_{\Omega} b_0  \ud x=0
\]
and
\[
 \frac{\ud}{\ud t}\int_{\Omega }\tilde a\ud x=0,
\]
which is the key to see that $\frac{1}{\eps^2}\int \tilde V_\eps^N\ud x \to 0$ when $\eps\to 0$.

The equation for the non-oscillating part can be decoupled from the rest of the system, and writes
\[
 \partial_t \Pi_0 \tilde V + \Pi_0 B(\Pi_0 \tilde V, \Pi_0\tilde V)=0,
\]
which can be rewritten, with the notation $\Pi_0\tilde V=(\bar a, \bar u, \bar b)$,
\[
 \partial_t \bar a + \bar u.\nabla_x \bar a=0,\qquad\qquad\partial_t \bar b+\bar u.\nabla_x \bar b=0,\qquad\qquad\nabla_x(b_0\bar a+\tau_0\bar b)=0,
\]
\[
 \partial_t \bar u + \bar u.\nabla_x \bar u + \nabla_x p=0,\qquad\nabla_x.\bar u=0.
\]
This is the incompressible Euler system. Now it remains to show that the approximate solution we constructed satisfies the constraint \eqref{constraint}. This is the object of the following result:

\begin{prop}
\label{A tend vers 0}
 With the above contruction of $\tilde V_\e^N =(a_1^N,u^N,b_1^N)$, the quantity
\be
A_\eps^N= \frac{1}{\eps^2}\frac{\ud}{\ud t}\iint_{\Omega\times\R^3}\frac{e^{-\eps b_1^N}}{b_0}|v|^2\frac{e^{a_0+\eps a_1^N-\frac{e^{-\eps b_1^N}|v|^2}{2b_0}}}{1+e^{a_0+\eps a_1^N-\frac{e^{-\eps b_1^N}|v|^2}{2b_0}}}\ud v\ud x
\ee
goes to 0 as $\eps\to 0$ for all $N\ge 1$.
\end{prop}

\begin{proof}
 $A_\eps^N$ can be rewritten
\[
 A_\eps^N=\frac{1}{\eps^2}\frac{\ud}{\ud t}\int (\beta_\eps^N)^{3/2}F_2(\alpha_\eps^N)\ud x 
\]
where
\[
 \beta_\eps^N = b_0 e^{\eps b_1^N},\qquad\qquad\alpha_\eps^N=e^{a_0+\eps a_1^N}
\]
and
\[
 F_p(\alpha)=\int |v|^p \frac{\alpha e^{-\frac{|v|^2}{2}}}{1+\alpha e^{-\frac{|v|^2}{2}}}\ud v.
\]
Using the equality
\[
 \alpha F_2'(\alpha)=\frac{3}{2}F_0(\alpha),
\]
we compute
\[
 A_\eps^N = \frac{3b_0^{3/2}}{2\eps}\int\left( e^{\frac{3}{2}\eps b_1^N}F_2(\alpha_\eps^N)\partial_t b_1^N + e^{\frac{3}{2}\eps b_1^N}F_0(\alpha_\eps^N)\partial_t a_1^N\right)\ud x.
\]
We easily prove the two identities
\[
 e^{\frac{3}{2}\eps b_1^N}(u.\nabla_x b_1^N + \frac{2}{3\eps}\nabla_x.u^N)= \nabla_x.\left(\frac{2}{3\eps}e^{\frac{3}{2}\eps b_1^N}u^N\right)
\]
and
\[
 \int F_2(\alpha_\eps^N)\nabla_x.\left(\frac{2}{3\eps}e^{\frac{3}{2}\eps b_1^N}u^N\right)\ud x = -\int e^{\frac{3}{2}\eps b_1^N}F_0(\alpha_\eps^N) u^N.\nabla_x a_1^N \ud x,
\]
which imply that
\[
\begin{split}
 A_\eps^N = &\frac{3b_0^{3/2}}{2\eps}\int \Big( F_2(\alpha_\eps^N)e^{\frac{3}{2}\eps b_1^N}\left(\partial_t b_1^N + u^N.\nabla_x b_1^N + \frac{2}{3\eps}\nabla_x.u^N\right) \\
&\qquad\qquad+ F_0(\alpha_\eps^N)e^{\frac{3}{2}\eps b_1^N}\left(\partial_t a_1^N + u^N.\nabla_x a_1^N\right)\Big)\ud x.
\end{split}
\]
Using the fact that $\mc A(a_1^N,u^N,b_1^N)\to 0$ in $L^2$ and the conservation of the mass, we deduce that 
\[
 A_\eps^N \underset{\eps \to 0}{\longrightarrow} 0.
\]

\end{proof}

\section{Useful intermediate results}
\label{intermediate}
In this section we define and study the linearized collision operator, and then give some bounds that will be very useful in the proof of
Theorem \ref{limite}.

\subsection{The linearized collision operator}

As we look at solutions of \eqref{scaled} which are fluctuations around an equilibrium state, that is,
\[
 f_\eps = P_\eps(1+\eps \bar g_\eps),
\]
it makes sense to use the linearized collision operator, which is defined as

\bee
\begin{split}
 \mc L_\eps g&=\int_{S^2}\int_{\R^3}B(v-v_*,\omega)\frac{M_{\eps,*}}{1+M_{\eps,*}}\frac{1}{1+M_\eps'}\frac{1}{1+M_{\eps,*}'}\\
&\quad\qquad\times\left(g(1+M_\eps)+g_*(1+M_{\eps,*})-g'(1+M_\eps')-g_*'(1+M_{\eps,*}')\right)\ud v_*\ud\omega.
\end{split}
\eee
Let us note
\be
\label{nu}
 \nu_\eps (v)=\int_{S^2}\int_{\R^3}\frac{B(v-v_*)M_{\eps,*}(1+M_\eps)}{(1+M_{\eps,*})(1+M_\eps')(1+M_{\eps,*}')}\ud v_*\ud\omega,
\ee
and recall that $P_\e =\frac{M_\e}{1+M_\e}$ with $M_\e=a_\e e^{-\frac{|v-u_\e|^2}{2b_\e}}$. Keep in mind that we will use $\mc L_\e$ with the collision kernel $B=\mathbbm{1}$.

%% ; the following proposition gives some properties
%% of $\mc L_\e$:

\begin{prop}
\label{linear}
Assume that the collision kernel $B(z,\omega)$ satisfies Grad's cutoff assumptions \eqref{Grad1}-\eqref{Grad2}; then $\mc L_\eps$ is a non-negative self-adjoint operator on $L^2(M_\eps\ud v)$ with domain
\[
 \mc D(\mc L_\eps) = \{g\in L^2(M_\eps \ud v)~|~\nu_\eps g \in L^2(M_\eps \ud v)\}=L^2(\nu_\eps M_\eps \ud v).
\]
and kernel
\[
 \textrm{Ker} \mc L_\eps = \textrm{span}\left\{\frac{1}{1+M_\eps},\frac{v}{1+M_\eps},\frac{|v|^2}{1+M_\eps}\right\}.
\]
It can be decomposed as
\[
 \mc L_\eps g(v) = \nu_\eps(v)g(v)+\mc K_\eps g(v)
\]
where $\mc K_\eps$ is a compact integral operator on $L^2(M_\eps \ud v)$ and
\[
 0<\nu_-\le\nu_\eps (v)\le \nu_+ (1+|v|^\beta)
\]
provided that 
\be
\label{abC}
\frac{1}{C_0}\le a_\eps \le C_0 ~~\textrm{and}~~ \frac{1}{C_0}\le b_\e \le C_0
\ee
 with $C_0$ a positive constant independent of $\e$, and $\nu_+, \nu_-$ depend only on $C_0$ and $C_B$ from \eqref{Grad2} (for $B=\mathbbm{1}$ we have $\nu_-=\frac{8\sqrt{2} \pi^{5/2}}{(1+C_0)^{9/2}}$ and $\nu_+=4\pi C_0^{9/2}$).
\end{prop}

\begin{rem}
In the framework of Theorem \ref{limite}, assumption \eqref{abC} is satisfied, since we consider for $a_\e$ and $b_\e$ small perturbations around the constant values $e^{a_0}$ and $b_0$.
\end{rem}

\begin{proof}
 Let $g,h \in L^2(M_\eps \ud v)$. Then
\[
\begin{split}
 \int h &\mc L_\eps g M_\eps \ud v= \int h(1+M_\eps) \mc L_\eps g \frac{M_\eps}{1+M_\eps}\ud v\\
&=\frac{1}{4}\iint\int_{S^2}B(v-v_*,\omega)\frac{M_\eps}{1+M_\eps}\frac{M_{\eps,*}}{1+M_{\eps,*}}\frac{1}{1+M_\eps'}\frac{1}{1+M_{\eps,*}'}\\
&\qquad\times\left(g(1+M_\eps)+g_*(1+M_{\eps,*})-g'(1+M_\eps')-g_*'(1+M_{\eps,*}')\right)\\
&\qquad\times\left(h(1+M_\eps)+h_*(1+M_{\eps,*})-h'(1+M_\eps')-h_*'(1+M_{\eps,*}')\right)\ud v\ud v_*\ud\omega
\end{split}
\]
which easily shows that $\mc L_\eps$ is a non-negative self-adjoint operator. Moreover, letting $h=g$ implies that $g$ is in the nullspace of $\mc L_\eps$ if and only if, for almost all $(v,v_*,\omega)\in \R^3\times\R^3\times S^2$,
\[
g(1+M_\eps)+g_*(1+M_{\eps,*})=g'(1+M_\eps')+g_*'(1+M_{\eps,*}').
\]
In other words, $g(1+M_\eps)$ must be a collision invariant, and therefore the nullspace of $\mc L_\eps$ is made of all linear
combinations of $\frac{1}{1+M_\eps}, \frac{v_i}{1+M_\eps}, \frac{|v|^2}{1+M_\eps}$ (see \cite{cercignani_fermions}).

The next step is to split $\mc L_\eps$ into two parts:
\[
 \mc L_\eps g=\nu_\eps (v)g +\mc K_\eps h(v)
\]
with $\nu_\e$ defined by \eqref{nu} and
\[
\begin{split}
 \mc K_\eps g(v)&=\int_{S^2}\int_{\R^3}B(v-v_*,\omega)\frac{M_{\eps,*}}{1+M_{\eps,*}}\frac{1}{1+M_\eps'}\frac{1}{1+M_{\eps,*}'}\\
&\qquad\qquad\qquad\times\left(g_*(1+M_{\eps,*})-g'(1+M_\eps')-g_*'(1+M_{\eps,*}')\right)\ud v_*\ud\omega.
\end{split}
\]
The bounds on $\nu_\eps$ can be proved using Grad's cutoff assumptions \eqref{Grad1}-\eqref{Grad2}. They are uniform in $\eps$ thanks to
\eqref{abC}. The operator $\mc K_\eps $ can be split into two operators in the following way:
\[
 \mc K_\eps g(v)= \mc K_\eps^1g(v) -\mc K_\eps^2 g(v)
\]
with
\[
 \mc K_\eps^1 g(v) = \int_{\R^3}\int_{S^2}B(v-v_*,\omega)M_{\eps,*}\frac{1}{1+M_\eps'}\frac{1}{1+M_{\eps,*}'}g_*\ud v_*\ud\omega
\]
and
\bee
\begin{split}
 &\mc K_\eps^2 g(v) =\int_{S^2}\int_{\R^3}B(v-v_*,\omega)\frac{M_{\eps,*}}{1+M_{\eps,*}}\frac{1}{1+M_\eps'}\frac{1}{1+M_{\eps,*}'}\\
&\qquad\qquad\qquad\times\left(g'(1+M_\eps')+g_*'(1+M_{\eps,*}')\right)\ud v_*\ud\omega.
\end{split}
\eee
In the classical case, the same decomposition holds, and a clever change of variables known as ``Carleman's parametrization'' allows to show that the following two operators are compact on $L^2(M_\eps \ud v)$ (see \cite{saintraybook} for example):
\[
 \bar {\mc K_1} g(v)= \int_{\R^3}\int_{S^2}B(v-v_*,\omega)M_{\eps,*}g_*\ud v_*\ud\omega
\]
\[
 \bar {\mc K_2} g(v)= \int_{\R^3}\int_{S^2}B(v-v_*,\omega)M_{\eps,*}(g'+g_*')\ud v_*\ud\omega.
\]
It is easy to see that
\[
 \|{\mc K_\eps^1} g\|_{L^2(M_\eps \ud v)}\le  \|\bar {\mc K_1} g\|_{L^2(M_\eps \ud v)}
\]
and
\[
 \|\mc K_\eps^2 g\|_{L^2(M_\eps \ud v)}\le \|\bar {\mc K_2} g\|_{L^2(M_\eps \ud v)}+\|\bar {\mc K_2} (M_\eps g)\|_{L^2(M_\eps \ud v)}
\]
from which we deduce that $\mc K_\eps$ is a compact operator on $L^2(M_\eps \ud v)$.
\end{proof}

With these results, we can assert that $\mc L_\eps$ is coercive, and therefore is a Fredholm operator:

\begin{cor}
\label{spectral gap}
 There exists $C_{\mc L}>0$ such that, for each $g\in \mc D(\mc L_\eps) \cap \left(\textrm{Ker} (\mc L_\eps)\right)^\bot$,
\[
 \int g \mc L_\eps g M_\eps \ud v \ge C_{\mc L} \|g\|^2_{L^2(\nu_\eps M_\eps \ud v)}.
\]
If \eqref{abC} holds, then $C_{\mc L}$ depends only on $C_0$ and $C_B$, but not on $\e$.
\end{cor}

\begin{proof}
 This result is a consequence of Proposition \ref{linear}. Indeed, the multiplication operator $g\mapsto \nu_\eps g$ is self-adjoint on $L^2(\nu_\eps M_\eps \ud v)$ and has continuous spectrum, namely the numerical range of $\nu_\eps$ $[\nu_\eps^-,+\infty[$ where
\[
 \nu_-\leq\nu_\eps^-=\inf_{v\in\R^3}\nu_\eps(v).
\]
Then, since $\mc K_\eps$ is self-adjoint and compact on $L^2(M_\eps \ud v)$, Weyl's theorem ensures that the spectrum of $\mc L_\eps$ is
made of $[\nu_\eps^-,+\infty[$ and of a sequence of eigenvalues on the interval $[0,\nu_\eps^-]$ with $\nu_\eps^-$ as unique possible
accumulation point. Consequently, there exists a smallest eigenvalue $\lambda_\eps^1$ which is bounded by below by $\lambda^1$ if
\eqref{abC} holds, and the following spectral gap inequality is
satisfied for each $g\in \mc D(\mc L_\eps) \cap \left(\textrm{Ker} (\mc L_\eps)\right)^\bot$:
\[
 \int g \mc L_\eps g M_\eps \ud v \ge C_{\mc L} \|g\|^2_{L^2( M_\eps \ud v)}.
\]
Recalling that
\[
 \int g \mc L_\eps g M_\eps \ud v= \int \nu_\eps g^2 M_\eps \ud v - \int g\mc K_\eps g M_\eps \ud v
\]
and using the continuity of $\mc K_\eps$, we get the inequality with the weighted norm as stated above.
\end{proof}

\subsection{Useful bounds}

This subsection lists some a priori bounds which are needed for the proof of convergence in the hydrodynamic limit. The first estimate comes
from Young's inequality:

\begin{prop}
\label{young}
 For $z>-1$ we define the function
\[
 h(z)=(1+z)\log (1+z)-z. 
\]
It satisfies
\[
 h(z) \ge 0 \quad \textrm{for}~ z>-1
\]
and
\[
 p|z|\le \lambda h^*(p) + \frac{1}{\lambda}h(z)
\]
where $h^*$ is the Legendre transform of $h$
\[
 h^*(p) = \max_{z>-1}(pz-h(z)) = e^p-p-1.
\]

\end{prop}

The study of the function $h$ is motivated by the relation
\[
H(f_\e|P_\e) = \iint_{\Omega\times\R^3}P_\eps h\left(\frac{f_\e-P_\e}{P_\e}\right)\ud v\ud x +  \iint_{\Omega\times\R^3} (1-P_\e)
h\left(\frac{P_\e-f_\e}{1-P_\e}\right)\ud v\ud x ,
\]

\begin{proof}
The first property is immediate from the definitions of $h$. The second one comes from Young's inequality
\[
 pz\le  h^*(p) + h(z),
\]
 supplemented with the two inequalities
\[
 h(|z|)\le h(z) \quad z>-1
\]
\[
 h^*(\lambda p)\le \lambda^2 h^*(p) \quad p\ge 0,~\lambda \in [0,1].
\]
\end{proof}
It will be useful to work in $L^2$ since it is the natural space for the study of $\mc L_\e$, and then we will use renormalized fluctuations
instead of the natural one defined as $f_\e = P_\e (1+\e \bar g_\e)$:

\begin{prop}
\label{estimates1}
Let us define the renormalized fluctuations 
\[
 \hat g_\eps =\frac{1}{\eps}\left(\sqrt{\frac{f_\eps}{P_\eps}}-1\right)
\]
and
\[
 \hat h_\e=\frac{1}{\e}\left(\sqrt{\frac{1-f_\e}{1-P_\e}}-1\right).
\]
Then
\be
\label{gcontrol}
 \iint_{\Omega\times\R^3} P_\eps \hat g_\eps^2 \ud v\ud x + \iint_{\Omega\times\R^3}(1-P_\e)\hat h_\e^2\ud v\ud x \leq \frac{1}{\eps^2}H(f_\eps|P_\eps).
\ee
Moreover, if \eqref{abC} holds holds,
\be
\label{estiminfty}
 \|\e \sqrt{P_\e}\hat g_\e\|_{L^\infty(\ud v\ud x)} +  \|\e\hat h_\e\|_{L^\infty(\ud v\ud x)} \le 4+\sqrt{C_0}.
\ee

\end{prop}

\begin{proof}
 Estimate \eqref{gcontrol} is a direct consequence of the definition of $\hat g_\eps$ and $\hat h_\e$ and the following inequality:
\[
 x\log\frac{x}{y}-x+y\ge  (\sqrt{x}-\sqrt{y})^2,\qquad\forall x,y>0.
\]
Recalling that $0\le f_\eps\le 1$ and that $P_\eps$ is bounded indipendently of $\e$ the second point comes by direct inspection.
\end{proof}
In the $L^2$ setting, the collision kernel needs also to be renormalized. We repeat here the d\'efinition of the entropy dissipation:
\bee
\begin{split}
 D(f_\e)=\frac{1}{4}\int_{\R^3\times S^2}B_\e(v-v_*,\omega)&\left(f_\e'f_{\e,*}'(1-f_\e)(1-f_{\e,*})-f_\e f_{\e,*}(1-f_\e')(1-f_{\e,*}')\right)\\
&\cdot\log \left(\frac{f_\e'f_{\e,*}'(1-f_\e)(1-f_{\e,*}}{f_\e f_{\e,*}(1-f_\e')(1-f_{\e,*}')}\right)\ud v_* \ud\omega.
\end{split}
\eee

\begin{prop}
\label{estimates2}
Define the renormalized collision kernel by
\[
\begin{split}
 &\hat q_\eps = \frac{1}{P_\eps}\iint B_\eps\Lambda_\eps^{1/2}\left(\sqrt{f_\eps'f_{\eps,*}'(1-f_\eps)(1-f_{\eps,*})}-\sqrt{f_\eps f_{\eps,*}(1-f_\eps')(1-f_{\eps,*}')}\right)\ud v_*\ud\omega,
\end{split}
\]
with
\[
 \Lambda_\eps=\frac{M_\eps M_{\eps,*}}{(1+M_\eps)(1+M_{\eps,*})(1+M_\eps')(1+M_{\eps,*}')}.
\]
 Then
\[
 \|\hat q_\eps\|_{L^2(\nu_\eps^{-1}P_\eps\ud v)}^2 \leq  D(f_\eps)
\]
where $\nu_\eps$ is the collision frequency defined by \eqref{nu}.
\end{prop}

Note that the scaling condition \eqref{scal} and the H theorem imply the following bound on the entropy dissipation
\[
\frac{1}{\e^{q+3}}\int_0^T \int_\Omega D(f_\e)\ud x\ud s \le C.
\]

\begin{proof}
 This estimate is obtained by the Cauchy-Schwarz inequality:
\[
 \hat q_\eps^2 \leq \frac{1}{P_\eps^2}\left(\iint B_\eps G_\e^2 \ud v_*\ud\omega\right)\left(\iint B_\eps\Lambda_\eps\ud v_*\ud\omega\right)
\]
with 
\[
 G_\e=\sqrt{f_\eps'f_{\eps,*}'(1-f_\eps)(1-f_{\eps,*})}-\sqrt{f_\eps f_{\eps,*}(1-f_\eps')(1-f_{\eps,*}')}.
\]
We easily see that
\[
 \iint B_\eps\Lambda_\eps\ud v_*\ud\omega =\nu_\e \frac{P_\e}{1+M_\e}.
\]
Next, using the classical inequality
\[
 (x-y)\log\frac{x}{y}\ge 4 (\sqrt{x}-\sqrt{y})^2,\qquad\forall x,y>0,
\]
we get the result.
\end{proof}
Mixing together the previous estimates, we can prove a relaxation result. Define $\Pi_\eps$ as the orthogonal projection in $L^2(P_\e \ud v)$ on $\textrm{Ker~}\mc L_\eps$. The following proposition will be useful in the control of the flux terms: coupled with the control on the tails \eqref{control} and an interpolation argument, it will allow the control of the third moment of $f_\e$.

\begin{prop}
\label{estimates3}
Define $\hat g_\e$ and $\hat h_\e$ as in Proposition \ref{estimates1} and assume \eqref{abC}. Then
\[
\begin{split}
&\left \|\left(\frac{\hat g_\eps-\hat h_\e}{1+M_\e}-\Pi_\eps \frac{\hat g_\eps-\hat h_\e}{1+M_\e}\right)\right\|_{L^2(P_\eps\ud v)}\\
 &=  O(\e)\Bigg(\|\hat g_\e\|_{L^2(P_\e)}^2+\|\hat h_\e\|_{L^2}^2 + \e\frac{1}{\e^2}\int P_\e h\left(\frac{f_\e-P_\e}{P_\e}\right)\ud v\Bigg)+o(\e).
\end{split}
\]

\end{prop}

\begin{proof}
Plugging the identities
\[
 f_\e=P_\e (1+\e \hat g_\e)^2,\qquad 1-f_\e =  (1-P_\e)(1+\e \hat h_\e)^2
\]
into
\[
 \sqrt{f'_\eps f'_{\eps,*}(1-f_\eps)(1-f_{\eps,*})}
\]
gives
\[
 \sqrt{f'_\eps f'_{\eps,*}(1-f_\eps)(1-f_{\eps,*})}=\sqrt{P_\e'P_{\e,*}'(1-P_\e)(1-P_{\e,*})}(1+\e\hat g_\e')(1+\e\hat g_{\e,*}')(1+\e\hat h)(1+\e\hat h_{\e,*})
\]
which leads to the key following identity:
\be
\label{keyident}
 \begin{split}
& \sqrt{f'_\eps f'_{\eps,*}(1-f_\eps)(1-f_{\eps,*})}-\sqrt{f_\eps f_{\eps,*}(1-f'_{\eps})(1-f'_{\eps,*})}\\
&= \eps\Lambda_\eps ^{1/2}\left(\hat g_\e'-\hat h_\e ' + \hat g_{\e,*}'-\hat h_{\e,*}' - (\hat g_\e-\hat h_\e)-(\hat g_{\e,*}-\hat h_{\e,*})\right)\\
&\quad+ \e^2 \Lambda_\e^{1/2} \Big(\hat g_\e' \hat g_{\e,*}'-\hat g_\e \hat g_{\e,*} + \hat h_\e \hat h_{\e,*} -\hat h_\e' \hat h_{\e,*}' + \hat g_\e' \hat h_\e -\hat g_\e \hat h_\e'\\
&\qquad+ \hat g_\e' \hat h_{\e,*} -\hat g_\e \hat h_{\e,*}'+ \hat g_{\e,*}' \hat h_\e -\hat g_{\e,*} \hat h_\e' + \hat g_{\e,*}' \hat h_{\e,*} -\hat g_{\e,*} \hat h_{\e,*}'\Big)\\
&\quad+ \e^3 \Lambda_\e^{1/2} \Big(\hat g_\e' \hat g_{\e,*}'(\hat h_\e +\hat h_{\e,*})-\hat g_\e \hat g_{\e,*}(\hat h_\e' +\hat h_{\e,*}') \\
&\qquad+ \hat h_\e \hat h_{\e,*}(\hat g_\e'+\hat g_{\e,*}') -\hat h_\e' \hat h_{\e,*}'(\hat g_\e + \hat g_{\e,*})\Big)\\
&\quad+ \e^4 \Lambda_\e^{1/2} \Big(\hat g_\e' \hat g_{\e,*}' \hat h_\e \hat h_{\e,*} - \hat g_\e \hat g_{\e,*} \hat h_\e' \hat h_{\e,*}'\Big).
 \end{split}
\ee
Raising it to the square and dividing by $\e^2$, we get, thanks to \eqref{estiminfty}:
\be
\label{controll}
\begin{split}
&\frac{1}{C}\int_{\R^3}\frac{\hat g_\eps-\hat h_\e}{1+M_\e}\mc L_\e \frac{\hat g_\eps-\hat h_\e}{1+M_\e} M_\e \ud v \\
&\le \frac{1}{\e^2}\iiint \left( \sqrt{f'_\eps f'_{\eps,*}(1-f_\eps)(1-f_{\eps,*})}-\sqrt{f_\eps f_{\eps,*}(1-f'_{\eps})(1-f'_{\eps,*})}\right)^2 \ud v\ud v_* \ud\omega\\
&+\e^2(\|\hat g_\e\|_{L^2(P_\e\ud v)}^4+\|\hat h_\e\|_{L^2(\ud v)}^4)
\end{split}
\ee
where $C$ depends on $C_0$. Note that here we used the fact that, for a given integrable function u, we have (see \cite{ADVW} for more details)
\[
 \int_{\R^3\times S^2}u(v')\ud v\ud\omega \le C \int_{\R^3}u(v)\ud v.
\]
The next step is to decompose $B= B_\eps + (B-B_\eps)$, so that
\[
\begin{split}
 &\frac{1}{\eps^2}\iiint B \left(\sqrt{f_\eps'f_{\eps,*}'(1-f_\eps)(1-f_{\eps,*})}- \sqrt{f_\eps f_{\eps,*}(1-f_{\eps}')(1-f_{\eps,*}')} \right)^2 \ud\omega\ud v\ud v_*\\
&\leq \frac{C}{\eps^2}D(f_\eps) + \frac{C}{\eps^2} \iiint (B-B_\eps)f_\eps f_{\eps,*} \ud\omega\ud v\ud v_*.
\end{split}
\]
Recalling that $f_\eps$ can be written $f_\eps=P_\eps(1+2 \eps \hat g_\eps + \eps^2 \hat g_\eps^2)$ and $M_\e=a_\e e^{-\frac{|v-u\e|^2}{2b_\e}}$, and thanks to the following inequality
\be
\label{petit}
 \mathbbm{1}_{|v-v_*|>\frac{1}{\eps^2}}M_\eps M_{\eps,*}\leq a_\eps^{1/2}M_\eps^{3/4}M_{\eps,*}^{3/4}e^{-\frac{1}{16 C_0\eps^4}},
\ee
we get
\be
\label{blabla}
\begin{split}
 &\frac{C}{\eps^2} \iiint \mathbbm{1}_{|v-v_*|>\frac{1}{\eps^2}}f_\eps f_{\eps,*} \ud\omega\ud v\ud v_*\\
&\le \frac{C}{\e^2}e^{-\frac{1}{C\e^4}} + 2C\iiint \mathbbm{1}_{|v-v_*|>\frac{1}{\eps^2}}P_\e P_{\e,*}\hat g_\e ^2\ud v\ud v_* \\
&\qquad\qquad +  4\e C\iiint \mathbbm{1}_{|v-v_*|>\frac{1}{\eps^2}}P_\e P_{\e,*}\hat g_\e ^2 \hat g_{\e,*}\ud v\ud v_* + C\e^2 \|\hat g_\e\|_{L^2(P_\e \ud v)}^4.
\end{split}
\ee
Using Proposition \ref{young} we obtain for $\delta >0$
\[
\begin{split}
 &\iiint \mathbbm{1}_{|v-v_*|>\frac{1}{\eps^2}}P_\e P_{\e,*}\hat g_\e ^2\ud v\ud v_* \\
&\le \delta\iint P_{\e,*} P_\e h\left(\frac{f_\e-P_\e}{P_\e}\right)\ud v\ud v_* + \iint \mathbbm{1}_{|v-v_*|>\frac{1}{\eps^2}}P_\e P_{\e,*} h^* \left( \frac{1}{\delta \e^2}\right)\ud v\ud v_*.
\end{split}
\]
We choose $\delta=20 C_0 \e^2$ :
\[
\begin{split}
 &\iiint \mathbbm{1}_{|v-v_*|>\frac{1}{\eps^2}}P_\e P_{\e,*}\hat g_\e ^2\ud v\ud v_* \\
&\le C \e^4 \left[ \left(\frac{1}{\e^2}\iint P_{\e,*} P_\e h\left(\frac{f_\e-P_\e}{P_\e}\right)\ud v\ud v_*\right)^2 \right]^{1/2} + o(1).
\end{split}
\]
We now apply the following variant of Young's inequality
\[
Y^\alpha \le \left( \frac{1}{\beta}\right)^{\frac{\alpha}{1-\alpha}} + \beta Y,\qquad\qquad Y\ge0, \beta >0, \alpha \in (0,1)
\]
with $\beta =1$, $\alpha=\frac{1}{2}$ and $Y=\left(\frac{1}{\e^2}\iint P_{\e,*} P_\e h\left(\frac{f_\e-P_\e}{P_\e}\right)\ud v\ud v_*\right)^2$:
\[
\begin{split}
 &\iiint \mathbbm{1}_{|v-v_*|>\frac{1}{\eps^2}}P_\e P_{\e,*}\hat g_\e ^2\ud v\ud v_* \\
&\le C \e^4  \left(\frac{1}{\e^2}\iint P_{\e,*} P_\e h\left(\frac{f_\e-P_\e}{P_\e}\right)\ud v\ud v_*\right)^2  + o(1).
\end{split}
\]

Since $\e \sqrt{P_{\e,*}}\hat g_{\e,*}$ is bounded, the same computation works for the second term in the right hand side of \eqref{blabla}, so that we finally get

\[
 \begin{split}
 &\frac{1}{\eps^2}\iiint B \left(\sqrt{f_\eps'f_{\eps,*}'(1-f_\eps)(1-f_{\eps,*})}- \sqrt{f_\eps f_{\eps,*}(1-f_{\eps}')(1-f_{\eps,*}')} \right)^2 \ud\omega\ud v\ud v_*\\
&\leq C \eps^2  \|\hat g_\eps\|_{L^2(P_\eps\ud v)}^4 + C\e^4 \left(\frac{1}{\e^2}\int P_\e h\left(\frac{f_\e-P_\e}{P_\e}\right)\ud v\right)^2 + o(1)
\end{split}
\]
since $\frac{1}{\e^2}D(f_\e)= o(1)$. Pluging this inequality into \eqref{controll} and using the coercivity of the linearized operator (Corollary \ref{spectral gap}) gives the result. 
\end{proof}

\section{End of the proof}
\label{sec:end}
We now have all the ingredients to prove Theorem \ref{limite}. We begin by controlling the flux terms, and then end the proof of the theorem.

\subsection{Control of the flux terms}

Our goal is to estimate the terms
\[
\frac{1}{\eps^2}\int_0^t\iint_{\Omega\times\R^3}f_\eps b_\eps\Psi_\eps\cdot\nabla_xb_1\ud x\ud v\ud s
\]
and
\[ 
 \frac{1}{\eps^2}\int_0^t\iint_{\Omega\times\R^3}f_\eps (\Phi_\eps\cdot\nabla_x u)\ud x\ud v\ud s
\]
with respect to the modulated entropy and the entropy dissipation. The main difficulty is that $\Psi_\eps = O(|v|^3)$ whereas the modulated entropy allows us to control the powers of $v$ only up to $2$, via Young's inequality. We therefore try to gain as much integrability as possible. Here the relaxation estimate \eqref{spectral gap} plays an central role. Since it gives a control on $\frac{\hat g_\eps -\hat h_\e}{1+M_\e}$, we use a new decomposition of $f_\e$

\be
\label{fg}
 f_\eps=P_\eps + \e^2 P_\e \left( \frac{\hat g_\e-\hat h_\e}{1+M_\e}\right)\left(\hat g_\e+\hat h_\e + \frac{2}{\e}\right)
\ee
and we obtain the following estimate:% thanks to the properties of the linearized collision operator:

\begin{prop}
\label{control flux}
 Assume that $\tilde \Phi_\eps$ and $\tilde \Psi_\eps$ have at most polynomial growth and that \eqref{abC} holds. Then the flux terms are bounded by the following quantities:
\bee
\begin{split}
& \frac{1}{\eps^2}\int_0^t\iint_{\Omega\times\R^3}f_\eps b_\eps\Psi_\eps.\nabla_xb_1\ud x\ud v\ud s +  \frac{1}{\eps^2}\int_0^t\iint_{\Omega\times\R^3}f_\eps (\Phi_\eps\cdot\nabla_x u)\ud x\ud v\ud s\\
&\le \frac{C}{\eps^2}\int_0^t \|D_x(u,b_1)(s)\|_{L^2\cap L^\infty(\Omega)}\left(H(f_\eps|P_\eps)(s)+\int_{\Omega}D(f_\eps)(s)\ud x\right)\ud s +o(1).
\end{split}
\eee
\end{prop}

\begin{proof}
The proof for the momentum flux (involving $\Phi_\eps$) is identical to the one with the energy flux (involving $\Psi_\eps$), hence we will focus on the first one. Writing $f_\eps$ as in \eqref{fg}, we use a first decomposition of the flux terms:
\be
\label{flux}
\frac{1}{\eps^2}\int_{\R^3} \Phi_\eps f_\eps \ud v = \frac{2}{\e}\int \Phi_\e P_\e \frac{\hat g_\e-\hat h_\e}{1+M_\e}\ud v + \int \Phi_\e P_\e \frac{\hat g_\e-\hat h_\e}{1+M_\e}(\hat g_\e +\hat h_\e)\ud v
\ee
since $\int_{\R^3}P_\eps \Phi_\eps\ud v=0$. To treat the first term in the right hand side, we use the identity \eqref{keyident}, wich implies
\be
\label{decompose}
\begin{split}
 &\frac{1}{\eps}\int_{\R^3}P_\eps \frac{\hat g_\e-\hat h_\e}{1+M_\e} \Phi_\eps \ud v = \int_{\R^3}M_\eps \tilde\Phi_\eps\frac{1}{\eps}\mc L_\eps\frac{\hat g_\e-\hat h_\e}{1+M_\e} \ud v\\
&=\frac{1}{\eps^2}\iiint B\tilde \Phi_\eps \Lambda_\eps^{1/2} \left(\sqrt{f_\eps'f_{\eps,*}'(1-f_\eps)(1-f_{\eps,*})}-\sqrt{f_\eps f_{\eps,*}(1-f_\eps')(1-f_{\eps,*}')}\right)\ud v\ud v_*\ud\omega\\
&- \iiint B\tilde \Phi_\eps\ \Lambda_\e \Big(\hat g_\e' \hat g_{\e,*}'-\hat g_\e \hat g_{\e,*} + \hat h_\e \hat h_{\e,*} -\hat h_\e' \hat h_{\e,*}' + \hat g_\e' \hat h_\e -\hat g_\e \hat h_\e'\\
&\qquad+ \hat g_\e' \hat h_{\e,*} -\hat g_\e \hat h_{\e,*}'+ \hat g_{\e,*}' \hat h_\e -\hat g_{\e,*} \hat h_\e' + \hat g_{\e,*}' \hat h_{\e,*} -\hat g_{\e,*} \hat h_{\e,*}'\Big)\ud v\ud v_*\ud\omega\\
&- \e \iiint B\tilde \Phi_\eps\Lambda_\e \Big(\hat g_\e' \hat g_{\e,*}'(\hat h_\e +\hat h_{\e,*})-\hat g_\e \hat g_{\e,*}(\hat h_\e' +\hat h_{\e,*}') \\
&\qquad+ \hat h_\e \hat h_{\e,*}(\hat g_\e'+\hat g_{\e,*}') -\hat h_\e' \hat h_{\e,*}'(\hat g_\e + \hat g_{\e,*})\Big)\ud v\ud v_*\ud\omega\\
&- \e^2 \iiint B\tilde \Phi_\eps\Lambda_\e \Big(\hat g_\e' \hat g_{\e,*}' \hat h_\e \hat h_{\e,*} - \hat g_\e \hat g_{\e,*} \hat h_\e' \hat h_{\e,*}'\Big)\ud v\ud v_*\ud\omega.
\end{split}
\ee
The first term can be estimated by the entropy dissipation as showed by the next lemma, while the other terms are controlled straightforwardly by the $L^2(M_\eps\ud v)$ norm of $\hat g_\eps$, and hence by the modulated entropy.

\begin{lem}
  Assume that $\tilde \Phi_\eps$ has at most a polynomial growth and that \eqref{abC} holds. Then the first term in the decomposition \eqref{decompose} is estimated by: 
%HERE WE NEED $\Omega$ BOUNDED
\bee
\begin{split}
& \int_\Omega \Bigg(\int B\tilde \Phi_\eps \Lambda_\eps^{1/2} \left(\sqrt{f_\eps'f_{\eps,*}'(1-f_\eps)(1-f_{\eps,*})}-\sqrt{f_\eps f_{\eps,*}(1-f_\eps')(1-f_{\eps,*}')}\right)\ud v\ud v_*\ud\omega\Bigg)^2\ud x\\
\leq& C\Bigg(\int_\Omega D(f_\eps)\ud x  + e^{-\frac{1}{16C_0\eps^4}} \Bigg)
\end{split}
\eee
\end{lem}

\begin{proof}

Let 
\[
 G_\e=\sqrt{f_\eps'f_{\eps,*}'(1-f_\eps)(1-f_{\eps,*})}-\sqrt{f_\eps f_{\eps,*}(1-f_\eps')(1-f_{\eps,*}')};
\]
We decompose
\be
\label{decompose2}
\begin{split}
 &\int B\tilde \Phi_\eps \Lambda_\eps^{1/2} G_\e\ud v\ud v_*\ud\omega\\
&=\int  B_\eps\tilde \Phi_\eps \Lambda_\eps^{1/2} G_\e\ud v\ud v_*\ud\omega+ \int  (B-B_\eps)\tilde \Phi_\eps \Lambda_\eps^{1/2} G_\e\ud v\ud v_*\ud\omega.
\end{split}
\ee
The first term is dealt with thanks to the Cauchy-Schwarz inequality and Proposition (\ref{estimates2}):
\[
\begin{split}
 &\int  B_\eps\tilde \Phi_\eps \Lambda_\eps^{1/2} G_\e\ud v\ud v_*\ud\omega\leq C\left(\int B (\tilde \Phi_\eps)^2 \Lambda_\eps \ud v_*\ud v\ud\omega\right)^{1/2} D(f_\eps)^{1/2}.
\end{split}
\]

By Corollary \ref{spectral gap}, the coercivity inequality implies, since $\tilde \Phi_\eps \in (\textrm{Ker}\mc L_\eps)^\bot$, that
\[
 \int B (\tilde \Phi_\eps)^2 \Lambda_\eps \ud v_*\ud v\ud\omega \le C \int \frac{M_\eps}{(1+M_\eps)^2}(\Phi_\eps)^2\ud v,
\]
with $C>0$ independent of $\eps$.

The second term in the right hand side of \eqref{decompose2}, containing the high velocities, is handled thanks to the inequality
\[
 \mathbbm{1}_{|v-v_*|>\frac{1}{\eps}}M_\eps M_{\eps,*}\leq a_\eps^{1/2}M_\eps^{3/4}M_{\eps,*}^{3/4}e^{-\frac{1}{16 b_\eps\eps^2}},
\]
using that $\tilde \Phi_\eps (1+M_\eps)$ is at most polynomial: 
\bee
\begin{split}
&\int  (B-B_\eps)\tilde \Phi_\eps (1+M_\eps)\Lambda_\eps^{1/2} G\ud v\ud v_*\ud\omega\\
&\leq Ce^{-\frac{1}{32 b_\eps\eps^4}}.
\end{split}
\eee

\end{proof}
The next terms in \eqref{decompose} are dealt with using Cauchy-Schwarz inequality and the a priori bound \eqref{estiminfty}, and therefore are easily bounded by $C(\|\hat g_\e\|_{L^2[P_\e\ud v}^2+\|\hat h_\e\|_{L^2}^2)$.
%\bee
%\begin{split}
% \int B& \tilde\Phi_\eps(1+M_\eps)\Lambda_\eps(\hat g_\eps \hat g_{\eps,*}-\hat g_\eps '\hat g_{\eps,*}')\ud v\ud v_*\ud\omega\ud x\\
%&\le \int_{\Omega}\left[\left\|\tilde \Phi_\eps\right\|_{L^2(M_\eps\ud v)}\left\|\iint B \frac{(1+M_\eps)}{M_\eps}\Lambda_\eps(\hat g_\eps \hat %g_{\eps,*}-\hat g_\eps'\hat g_{\eps,*}')\ud v_*\ud\omega\right\|_{L^2(M_\eps\ud v)}\right]\ud x\\
%&\le C\int_{\Omega}\|\Phi_\eps\|_{L^2(M_\eps\ud v)}\|\hat g_\eps\|_{L^2(M_\eps\ud v)}^2\ud x\\
%&\le \frac{C}{\eps^2}H(f_\eps|P_\eps)
%\end{split}
%\eee
%
Coming back to \eqref{flux}, the second term on the right hand side is splitted as follows, in order to use the relaxation estimate:
\[
\begin{split}
 \int \Phi_\e P_\e \frac{\hat g_\e-\hat h_\e}{1+M_\e}(\hat g_\e +\hat h_\e)\ud v=& \int \Phi_\e P_\e \left(\frac{\hat g_\e-\hat h_\e}{1+M_\e}-\Pi_\e\frac{\hat g_\e-\hat h_\e}{1+M_\e}\right)(\hat g_\e +\hat h_\e)\ud v\\
&\qquad+\int \Phi_\e P_\e \Pi_\e\left(\frac{\hat g_\e-\hat h_\e}{1+M_\e}\right)(\hat g_\e +\hat h_\e)\ud v
\end{split}
\]
The last term is easily handled:
\[
 \int \Phi_\e P_\e \Pi_\e\left(\frac{\hat g_\e-\hat h_\e}{1+M_\e}\right)(\hat g_\e +\hat h_\e)\ud v \leq C(\|\hat g_\e\|_{L^2(P_\e \ud v)} + \|\hat h_\e\|_{L^2})
\]
while we need some control on the high velocities for the first one. From the hypothesis \eqref{control}, we deduce that
\[
 \int_{\R^3}P_\e \left(\frac{f_\e}{P_\e}\right)^{2p}\ud v \leq C_p~~\textrm{a.e.}
\]
for all $p<1$ and uniformly in $\e$ (which depends on $p$), since the moments of $P_\e$ differ from that of $P_0$ only by terms of order $\e$. This estimate leads to some control on the rescaled fluctuations:
\[
 \e |\hat g_\e| +\e |\hat h_\e|=O(1)_{L^\infty_{t,x}L^{4p}(P_\e\ud v)}
\]
so that
\be
\label{estim4}
 \e\left(\left(\frac{\hat g_\e-\hat h_\e}{1+M_\e}\right)-\Pi_\e\frac{\hat g_\e-\hat h_\e}{1+M_\e}\right)=O(1)_{L^\infty_{t,x}L^{4p}(P_\e\ud v)}.
\ee
%The relaxation estimate \eqref{spectral gap} coupled with some interpolation gives
We easily find the following interpolation inequality for all functions $\phi \in L^2 \cap L^{4p}$
\[
\|\phi\|_{L^{\frac{8p}{2p+1}}}\le \|\phi\|_{L^2}^{1/2 }\|\phi\|_{L^{4p}}^{1/2},
\]
which implies, together with  \eqref{estim4} and Proposition \ref{estimates3}:
\[
\begin{split}
&\left\|\frac{\hat g_\e-\hat h_\e}{1+M_\e} -\Pi_\e\frac{\hat g_\e-\hat h_\e}{1+M_\e} \right\|_{L^{\frac{8p}{2p+1}}(P_\e \ud v)}\\
&= O\bigg(\|\hat g_\e\|_{L^2(P_\e \ud v)} + \|\hat h_\e\|_{L^2} + \left( \frac{1}{\e^2} \int P_\e h\left(\frac{f_\e-P_\e}{P_\e}\right)\ud v\right)^{1/2}\bigg) +o(1).
\end{split}
\]
Choosing $p>\frac{1}{2}$ implies $\frac{8p}{2p+1}>2$, and H\"older's inequality gives
\[
\int \phi_1 \phi_2 \phi_3 P_\e \ud v \le \|\phi_1\|_{L^2(P_\e)} \|\phi_2\|_{L^{\frac{8p}{2p+1}}(P_\e)} \|\phi_3\|_{\frac{8p}{2p-1}(P_\e)}
\]
such that taking $\phi_1=\hat g_\e +\hat h_\e$, $\phi_2=\frac{\hat g_\e-\hat h_\e}{1+M_\e} -\Pi_\e\frac{\hat g_\e-\hat h_\e}{1+M_\e}$ and $\phi_3=\Phi_\e$ leads to the final estimate:
\[
\begin{split}
 \int \Phi_\e P_\e \left(\frac{\hat g_\e-\hat h_\e}{1+M_\e}-\Pi_\e\frac{\hat g_\e-\hat h_\e}{1+M_\e}\right)(\hat g_\e +\hat h_\e)\ud v&= O\bigg(\|\hat g_\e\|_{L^2(P_\e \ud v)}^2 + \|\hat h_\e\|_{L^2}^2\\
&\qquad+  \frac{1}{\e^2} \int P_\e h\left(\frac{f_\e-P_\e}{P_\e}\right)\ud v\bigg) +o(1) 
\end{split}
\]

\end{proof}

\subsection{Proof of convergence}

In view of the preceding results, we are now able to prove Theorem \ref{limite}. Proposition \ref{modulated entropy identity} together with Proposition \ref{control flux} imply the following:
\[
 \begin{split}
&\frac{1}{\eps^2}H(f_\eps|P_\eps)(t)+\frac{1}{\eps^{q+3}}\int_0^t\int_{\Omega}D(f_\eps)\ud x\ud s\\
&\le \frac{1}{\eps^2}H(f_{\eps,in}|P_{\eps,in}) -\frac{1}{\eps}\int_0^t\iint_{\Omega\times\R^3}f_\eps\left(1,\frac{v-\eps u}{b_\eps}, \frac{|v-\eps u|^2}{2b_\eps}\right).\mc A_\eps(a_1,u,b_1)\ud x\ud v \ud s\\
&+\frac{1}{3\eps^2}\int_0^t\left[\frac{\ud}{\ud t}\iint_{\Omega\times\R^3}\frac{e^{-\eps b_1}}{b_0}|v|^2\frac{e^{a_0+\eps a_1-\frac{e^{-\eps b_1}|v|^2}{2b_0}}}{1+e^{a_0+\eps a_1-\frac{e^{-\eps b_1}|v|^2}{2b_0}}}\ud v\ud x \right]\ud s\\
&+\frac{C}{\eps^2}\int_0^t \|D_x(u,b_1)(s)\|_{L^2\cap L^\infty(\Omega)}\left(H(f_\eps|P_\eps)(s)+\int_{\Omega}D(f_\eps)(s)\ud x\right)\ud s +o(1).\\
\end{split}
\]
We plug the functions $(a_1^{\eps,N}, u^{\eps,N}, b_1^{\eps,N})$ constructed in Theorem \ref{approxsols} in the approximate solution $P_\eps$. It follows by Gronwall's lemma that 
\[
\begin{split}
 \frac{1}{\eps^2}H&(f_\eps|P_\eps^N)(t) \leq  \frac{1}{\eps^2}H(f_{\e,in}|P_{\eps,in}^N)\exp\left(\int_0^t \|D_x(u^{\eps,N},b_1^{\eps,N})(s)\|_{L^2\cap L^\infty (\Omega)}\right)\ud s\\
&+\int_0^t \exp\left(\int_s^t \|D_x(u^{\eps,N},b_1^{\eps,N})(\tau)\|_{L^2\cap L^\infty (\Omega)}\ud \tau\right)\\
&\qquad\times \Bigg[\frac{1}{\eps}\iint_{\Omega\times\R^3}f_\eps\left(1,\frac{v-\eps u^{\eps,N}}{b_\eps}, \frac{|v-\eps u^{\eps,N}|^2}{2b_\eps}\right).\mc A_\eps(a_1^{\eps,N},u^{\eps,N},b_1^{\eps,N})\ud x\ud v \\
&\qquad\qquad\quad+\frac{1}{3\eps^2}\Bigg(\frac{\ud}{\ud t}\iint_{\Omega\times\R^3}\frac{e^{-\eps b_1^{\eps,N}}}{b_0}|v|^2\frac{e^{a_0+\eps a_1^{\eps,N}-\frac{e^{-\eps b_1^{\eps,N}}|v|^2}{2b_0}}}{1+e^{a_0+\eps a_1^{\eps,N}-\frac{e^{-\eps b_1^{\eps,N}}|v|^2}{2b_0}}}\ud v\ud x \Bigg) +o(1) \Bigg]\ud s. 
\end{split}
\]
The last term go to 0 as $\eps\to 0$ and $N\to 0$ thanks to Theorem \ref{approxsols}, using the fact that 
\[
 \int_\Omega \mc A_\eps(a_1^{\eps,N},u^{\eps,N},b_1^{\eps,N})\ud x=0.
\]
The fist term can be decomposed as
\[
\begin{split}
 \frac{1}{\eps^2}H(f_{\eps,in}|P_{\eps,in}^N)=& \frac{1}{\eps^2}H(f_{\eps,in}|P_{(a_1^{in},u^{in},b_1^{in})}) + \frac{1}{\eps^2}H(P_{(a_1^{in},u^{in},b_1^{in})}|P_{\eps,in}^N)\\
&+\frac{1}{\eps^2} \iint\Big( (f_{\eps,in}-P_{(a_1^{in},u^{in},b_1^{in})})\log \frac{P_{(a_1^{in},u^{in},b_1^{in})}}{P_{\eps,in}^N}\\
& \qquad+ (P_{(a_1^{in},u^{in},b_1^{in})}-f_{\eps,in})\log \frac{1-P_{(a_1^{in},u^{in},b_1^{in})}}{1-P_{\eps,in}^N}\Big) \ud v\ud x
\end{split}
\]
which tends to 0 as $\eps\to 0$ and $N\to\infty$ since $\frac{1}{\eps^2}H(f_{\eps,in}|P_{(a_1^{in},u^{in},b_1^{in})})\to 0$ by hypothesis and $(a_1^{\eps,N,in}, u^{\eps,N,in}, b_1^{\eps,N,in})\to (a_1^{in}, u^{in}, b_1^{in})$ by Theorem \ref{approxsols}.

This leads to
\[
 \frac{1}{\e^2}H(f_\e|P_\e^N)\to 0 ~~\textrm{in}~L^\infty_{loc}([0;T))~\textrm{as}~\e\to 0~\textrm{and}~N\to\infty,
\]
which is a stronger result respect to the one stated above.

\section*{Acknowledgment}
This work was partialy supported by the GDRE 224 CNRS-INdAM GREFI-MEFI and by Projet CBDif-Fr ANR-08-BLAN-0333-01. I thank the anonymous referee for his comments that helped a lot improving this manuscript.

\bibliographystyle{abbrv}
\bibliography{biblio}

\end{document}